\documentclass[11pt,eqsection]{article}

\usepackage{epsfig}
\usepackage{amsfonts,layout}
\usepackage{amsfonts,layout}
\begin{document}

\newtheorem{Theorem}{Theorem}[section]
\newtheorem{Definition}[Theorem]{Definition}
\newtheorem{Proposition}[Theorem]{Proposition}
\newtheorem{Lemma}[Theorem]{Lemma}
\newtheorem{Corollary}[Theorem]{Corollary}

\newtheorem{Remark}{Remark}[section]
\newtheorem{Example}{Example}[section]

\newcommand{\R}{I\hspace{-1.5 mm}R}
\newcommand{\N}{\mathbb N}
\newcommand{\Z}{\mathbb Z}
\def\K1{{\cal K}_1}
\def\Kes{{\cal K}_1^{\epsilon,\sigma}}
\def\HK2{\widehat{{\cal K}_2}}

\def\div{{\rm div\,}}

\def\QED{\begin{flushright}QED  \end{flushright}}

\title{Some flows in shape optimization}

\author{Pierre Cardaliaguet\thanks{Universit\'e de Bretagne
Occidentale, UFR des Sciences et Techniques, 6 Av. Le Gorgeu,
BP 809, 29285 Brest, France; e-mail:
$<${\tt Pierre.Cardaliaguet@univ-brest.fr}$>$ }
and Olivier Ley\thanks{Laboratoire de Math\'ematiques et Physique Th\'eorique.
Facult\'e des Sciences et Techniques, Parc de Grandmont, 37200 Tours, France;
e-mail: $<${\tt ley@gargan.math.univ-tours.fr}$>$}}

\date{}
\maketitle 

{\small
\noindent {\bf Abstract.}  Geometric flows related to 
shape optimization problems of Bernoulli type are investigated. The evolution law
is the sum of a curvature term and a nonlocal term of Hele-Shaw type. We introduce
generalized set solutions, the definition of which is widely inspired by
viscosity solutions. The main result is an inclusion preservation principle for 
generalized solutions. As a consequence, we obtain existence, uniqueness and stability
of solutions. Asymptotic behavior for the flow
is discussed: we prove that the solutions converge to a generalized Bernoulli exterior 
free boundary problem.
\vspace{3mm}

\noindent{\bf R\'esum\'e. } 
On \'etudie des flots g\'eom\'etriques li\'es \`a des probl\`emes d'optimisation
de forme du type ``probl\`{e}me de Bernoulli". La loi d'\'evolution consid\'er\'ee
a la forme d'une somme d'un terme de courbure et d'un terme non-local de type
Hele-Shaw. Notre d\'efinition de solution g\'en\'eralis\'ee est fortement inspir\'ee
de la notion de solutions de viscosit\'e. Le r\'esultat central est un principe
d'inclusion pour les ensembles solutions. Nous en d\'eduisons l'existence, une
unicit\'e g\'en\'erique et des propri\'et\'es de stabilit\'e des solutions.
Enfin, nous \'etudions le comportement asymptotique des solutions en montrant
qu'elles convergent vers la solution d'un probl\`eme \`a fronti\`ere libre
de Bernoulli.}

\section{Introduction}

In recent years several works have been devoted to the study of viscosity solution
for moving boundary problems whose evolution law is governed by a nonlocal equation.
See in particular \cite{af02, cr04, cardaliaguet00,cardaliaguet01, dlks04, kim03, kim04}.
In this paper, we consider 
subsets $\Omega(t)$ of $\R^N$ (with $N\geq2$) whose boundary $\partial \Omega(t)$ evolves with a normal
velocity of the type
\begin{equation}\label{veloci}
V_{(t,x)}^\Omega= F(\nu_x^{\Omega(t)}, H_x^{\Omega(t)})+ \lambda \bar{h}(x,\Omega(t))
\end{equation}
where $\lambda\geq0$, $\nu_x^{\Omega(t)}$ is the outward unit normal to $\partial \Omega(t)$
at $x$, $H_x^{\Omega(t)}$ is the curvature matrix of $\partial \Omega(t)$ at $x$ (nonpositive 
for convex sets), $F$ is continuous and elliptic, i.e., nondecreasing with respect to the curvature matrix. 
The nonlocal term $\bar{h}$ is of Hele-Shaw type: 
$$
\bar{h}(x,\Omega(t))=  |D  u(x)|^2\;,
$$
where $u:\Omega(t)\to \R$ is the solution to the following  p.d.e.
\begin{equation}\label{ppde}
\left\{\begin{array}{cl}
 -\Delta  u=0 & {\rm in }\; \Omega(t)\backslash S,\\
 u= g & {\rm on }\; \partial S,\\
 u=0 & {\rm on }\; \partial \Omega(t).
\end{array}\right.
\end{equation}
The set $S\neq\emptyset$ is a fixed source with a smooth boundary
and $g:\partial S\to \R$ is positive and smooth. We always assume 
that $S\subset\subset \Omega(t)$. 

The motivation to study such problems comes from several numerical works using the ``level-set
approach" in shape optimization \cite{ajt04, os01, sethian99, sw00, wwg03}. 
The idea of these papers is to use formally a gradient method for the minimization of 
an objective function $J(\Omega)$ where $\Omega$ is a subset of $\R^N$. 
The use of the level-set method for building the gradient flow has then the major
advantage to allow topological changes. Let us underline that this technique is up to now 
purely heuristic. One of the goals of this paper is to justify it for some simple
shape optimisation problems.

In order to make our purpose more transparent, a brief description of the level-set approach 
to shape optimization problem is now in order
(see also the discussion in \cite{ajt04} for a more detailed presentation concerning more
realistic shape optimization problems). Consider the problem
of minimizing the capacity of a set under volume constraints:
\begin{equation}\label{OptiPb1}
\min_{S\subset\subset \Omega\subset\subset \R^N}\left\{{\rm cap}(\Omega)\quad{\rm with}\quad 
{\rm vol}(\Omega)=\;{\rm constant}\right\}\;,
\end{equation}
where 
$$
{\rm cap}(\Omega)=\int_{\Omega\backslash S}|Du(x)|^2dx\quad{\rm and }\quad 
{\rm vol}(\Omega)=\int_{\Omega\backslash S}dx
$$
and $u$ is the solution of (\ref{ppde}) with $\Omega$ instead of $\Omega(t)$.
For any local diffeomorphism $\theta$, we can compute the shape derivatives with respect
to $\theta$ of the capacity and of the volume. By Hadamard formulas we get
$$
{\rm cap}'(\Omega)(\theta)=
\int_{\partial \Omega}|Du(\sigma)|^2\langle\theta(\sigma),\nu_\sigma^\Omega\rangle d\sigma
\ {\rm and }\ 
{\rm vol}'(\Omega)(\theta)=\int_{\partial \Omega}\langle\theta(\sigma),\nu_\sigma^\Omega\rangle d\sigma\;.
$$
Assuming that the optimal shape $\Omega$ is smooth, the necessary conditions of optimality
states that there is a Lagrange multiplier $\Lambda>0$ such that 
$$
{\rm cap}'(\Omega)(\theta)+\Lambda\, {\rm vol}'(\Omega)(\theta)=0\;.
$$
So it is natural to set
$$
J_\lambda (\Omega)={\rm vol}(\Omega)+\lambda\, {\rm cap}(\Omega)\;,
$$
where $\lambda=1/\Lambda$. If we choose $\theta(x)=(-1+\lambda |Du(x)|^2)\nu_x^\Omega$ on $\partial \Omega$, 
then, at least formally, we get 
$$
J'_\lambda(\Omega)(\theta)=-\int_{\partial \Omega}(-1+\lambda|Du(\sigma)|^2)^2d\sigma\leq 0\;.
$$
Therefore the velocity $\theta(x)=(-1+\lambda|Du(x)|^2)\nu_x^\Omega$ appears as a descent
direction for the optimization problem (\ref{OptiPb1}) and for the set $\Omega$. The heuristic method
for solving (\ref{OptiPb1}) is now clear: fix an initial position $\Omega_0$,
consider the evolution $(\Omega(t))_{t\geq0}$ with normal velocity given by (\ref{veloci})
and $F\equiv-1$, and compute the limit of $\Omega(t)$ as $t\to+\infty$: this limit
is the natural candidate minimizer for  (\ref{OptiPb1}).

It is worth noticing that problem (\ref{OptiPb1}) has for necessary condition the classical Bernoulli exterior
free boundary problem
\begin{eqnarray}\label{bernoulliEFB}
{\rm Find \ a \ set} \ K\subset\subset \R^N, \ {\rm with} \  S\subset\subset K \ 
{\rm and} \ |Du(x)|= k \ {\rm for \ all} \  x\in \partial K,
\end{eqnarray}
where $k>0$ is a fixed constant and $u$ is the solution of (\ref{ppde}). We refer the
reader to the survey paper \cite{fr97} for a complete description of this problem.

If one considers a perimeter constraint instead of a volume constraint:
$$
\min_{S\subset\subset \Omega\subset\subset \R^N}\left\{{\rm cap}(\Omega)\quad{\rm with}\quad 
{\rm per}(\Omega)=\;{\rm constant}\right\}\;,
$$
one is naturally lead to consider the evolution equation (\ref{veloci}) with 
$F(\nu,A)=\frac{1}{N-1}Tr(A)$ (i.e., the mean curvature). The flow is then formally a descent direction for 
$$
J_\lambda(\Omega)={\rm per}(\Omega)+\lambda {\rm cap}(\Omega)\;.
$$
Let us underline that this problem has for necessary condition the generalization of the free boundary
problem (\ref{bernoulliEFB}) with curvature dependance (see (\ref{pb-bernoulli})).

Of course all the above computations are only formal: in general, solutions to the evolution
equation do not remain smooth, even when starting from smooth initial data. Numerically, this difficulty 
is overcome by using the level-set approach, which allows to define the solution after the onset of singularities.
The aim of this paper is to define and study generalized solutions of the evolution equation, 
and to investigate the asymptotic behavior of the solution as $t\to+\infty$. 

Our concept of solutions is widely inspired by the definition of viscosity solution for the mean curvature motion,
which corresponds to equation (\ref{veloci}) with $F(\nu,A)=\frac{1}{N-1}Tr(A)$ and $\lambda=0$.
Motivated by the numerical work of Osher and Sethian \cite{os88},  
a weak notion of solution for this motion was introduced in the articles
of Chen, Giga and Goto \cite{cgg91} and Evans and Spruck \cite{es91}. 
In this so-called level-set method, 
the evolution is described as the level set of the solution of an auxiliary pde, the level set equation.
This equation is solved in the sense of viscosity solutions (see \cite{cil92}). This powerful
method leads to plenty of results, we refer for instance to the survey book of Giga \cite{giga02}.
Note that the level-set approach in shape optimization is a natural--but up to now formal--generalization 
of these ideas.

As pointed out in \cite{bss93, bs98, soner93}, the generalized solutions obtained by 
the level set approach can also be defined in more geometric and intric ways 
(see also the related notion of barrier solutions introduced by De Giorgi). 
We use here a definition introduced in \cite{af02}, and used repetitively
in \cite{cr04, cardaliaguet00, cardaliaguet01}. In the case of the mean curvature motion,
Giga \cite{giga02} proved this definition is equivalent to the level-set one. 
Compared with the already quoted studies on viscosity solutions of front
propagation problems with nonlocal terms, the main novelty of this paper is 
the fact that we are able to treat signed velocities which also involve curvature terms.
We learnt recently that a similar result (for a Stefan problem) has  been obtained by Kim 
in \cite{kim04}.

Our main result is an inclusion principle, which is the equivalent of the maximum
principle for geometric evolutions. It states that viscosity subsolutions for the flow remain included
into viscosity supersolutions, provided the initial positions are. For this we have to generalize
Ilmanen interposition Lemma, which was already the key tool of  \cite{ cr04, cardaliaguet01}.
This Lemma allows to separate disjoint sets by a smooth (that is ${\cal C}^{1,1}$)
surface in a clever way. We improve this result in two directions (see Theorem \ref{interposition}). 
At first we show that, when dealing with
subsets of $\R\times \R^N$, the smooth separating hypersurfaces in $\R\times \R^{N}$ 
can be chosen to be smoothly evolving hypersurfaces of $\R^N$. Secondly, 
we build in a carefull way a ${\cal C}^2$ approximation of  these evolving hypersurfaces which allows
to treat problems with curvature as in (\ref{veloci}).

Let us finally explain how this paper is organized. In Section 2, we define the notion
of generalized solutions and state the main properties of the velocity law.
Section 3 is devoted to the interposition Theorems. In Section 4 we state and prove the
inclusion principle for our generalized solutions. As a consequence, we derive results
about existence, uniqueness and stability of generalized solutions.
Finally, Section 5 is devoted to the asymptotic behaviour of the solutions
in terms of a generalized Bernoulli exterior free boundary problem. \\

\noindent{\bf Acknowledgment.} The authors are partially supported by the ACI grant
JC 1041 ``Mouvements d'interface avec termes non-locaux'' from the French Ministry
of Research.

\section{Definitions and preliminary results}

\subsection{Definition of the solutions}

Let us first fix some notations: throughout the paper $|\cdot|$ denotes
the euclidean norm (of $\R^N$ or $\R^{N+1}$, depending on the context)
and $B(x,R)$ the open ball centered at $x$ and of radius $R$.
If $K$ is a subset of $\R^N$ and $x\in\R^N$, then $d_K(x)$ denotes the usual
distance from $x$ to $K$: $d_K(x)=\inf_{y\in K} |y-x|$ and ${\bf d}_K$
is the signed distance to $\partial K$ defined by
\begin{eqnarray} \label{def-dist-signee}
{\bf d}_K (x)=\left\{\begin{array}{ll}
d_K (x) & {\rm if }\quad x\notin K, \\
-d_{\partial K}(x) & {\rm if }\quad x\in K.
\end{array}\right.
\end{eqnarray}
Finally, in the whole paper, if $K_1$ and $K_2$ are subset of $\R^M$ for $N\geq 1,$
then
\begin{eqnarray*}
K_1 \subset\subset K_2
\end{eqnarray*}
means that $K_1$ is bounded and, either $\overline{K_1}\subset {\rm int}(K_2)$ or equivalently 
$\overline{K_1}\cap \overline{\R^N\backslash K_2}=\emptyset.$

We intend to study the evolution of compact hypersurfaces $\Sigma(t)=\partial \Omega(t)$ of
$\R^N$, where $\Omega(t)$
is an open set, evolving with 
the following law:\begin{equation}\label{fpp}
\forall t\geq0 \; , \; x\in \Sigma(t), \; \; V_{(t,x)}^\Omega= h_\lambda(x, \Omega(t))
\end{equation}
where $V_{(t,x)}^\Omega$ is the normal velocity
of the evolving set, $h_\lambda=h_\lambda(x,\Omega)$ is  given, for any set $\Omega\subset\R^N$ with smooth boundary by
\begin{equation}\label{defh}
h_\lambda (x,\Omega)= F(\nu_x^\Omega, H_x^\Omega)+ \lambda \bar{h}(x,\Omega)
\end{equation}
where $\nu_x^\Omega$ is the outward unit normal to $\Omega$ at $x$, $H_x^\Omega$ the curvature
matrix. Throughout this paper we assume that $(\nu , A) \in S^{N-1}\times {\cal S}_N \mapsto F(\nu ,A)\in \R$
is continuous and elliptic, i.e., nondecreasing with respect to
the matrix. Here $S^{N-1}$ denotes the $(N-1)-$dimensional 
unit sphere, and ${\cal S}_N$ the space of $N-$dimensional symmetric matrices.
Typical examples for $F$ are $F (\nu , A)=-1$ (this corresponds to the flow associated to
Bernoulli problem in the introduction) or $F(\nu , A)={\rm Tr} (A)$ (for
the flow arising in the  minimization of the capacity under perimeter constraints).
As for $\bar{h}$, it is a nonlocal evolution term of Hele-Shaw type:
the example we consider here is
\begin{equation}\label{Defh2}
\bar{h}(x,\Omega)=  |D  u(x)|^2\;,
\end{equation}
where $u:\Omega\to \R$ is the solution of the following  p.d.e.
\begin{equation}\label{pde}
\left\{\begin{array}{lcl}
i) & -\Delta  u=0 & {\rm in }\; \Omega\backslash S,\\
ii) & u= g & {\rm on }\; \partial S,\\
iii) & u=0 & {\rm on }\; \partial \Omega.
\end{array}\right.
\end{equation}
The set $S\neq\emptyset$ is a fixed source and we always assume above that $S\subset\subset \Omega(t)$.
Here and throughout the paper, we suppose that 
\begin{equation}\label{Hypfg}
\left\{\begin{array}{rl}
i) & \mbox{\rm $S\subset \R^N$ is bounded and equal
to the closure of an open set}\\
& \mbox{\rm  with a ${\cal C}^2$ boundary,  }\\
ii) & \mbox{\rm $g:\partial S\to(0,+\infty)$ is ${\cal C}^{1,\alpha}$     (for some $\alpha\in(0,1)$). }
\end{array}\right.
\end{equation}
Let us underline that 
$\bar{h}(x,\Omega)$ is well defined as soon as $\Omega$ has a ``smooth" 
(say for instance ${\cal C}^{1,\alpha}$) boundary 
and that $S\subset\subset \Omega$. 
In the sequel, we set
\begin{equation}\label{DefD}
{\cal D}= \{ K\subset \R^N \; : \;K\ {\rm is \ bounded}\; {\rm and }\;  S\subset {\rm int}(K)\}\;,
\end{equation}
where ${\rm int}(K)$ denotes the interior of $K$.

From now on, we consider the graph 
$$
{\cal K}=\{(t,x)\in \R^+\times \R^N \ : \  x\in\Omega(t)\}\;.
$$ 
of the evolving sets $\Omega(t).$ Note that ${\cal K}$ is a subset of $\R^+\times \R^N$. 
The set ${\cal K}$ is our
main unknown. We denote by $(t,x)$ an
element of such a set, where $t\in\R^+$ denotes the time and $x\in\R^N$
denotes the space. We  set 
$$ 
{\cal K}(t) \;=\; \{x\in\R^N\; :\; (t,x)\in {\cal K}\}\;.
$$
The closure of the
set ${\cal K}$ in $\R^{N+1}$ is denoted by $\overline{\cal K}$.
The closure of the complementary of ${\cal K}$ is denoted $\widehat{{\cal K}}$:
$$
\widehat{{\cal K}}=\overline{\left(\R^+\times\R^N\right)\backslash {\cal K}}
$$
and we set
$$
\widehat{{\cal K}}(t)=\{x\in\R^N\; : \; (t,x)\in \widehat{{\cal K}} \}\;.
$$

We use here repetitively the terminology and the notations introduced in \cite{cardaliaguet00, cardaliaguet01} and  
\cite{cr04}: 
\begin{itemize}

\item
A {\it tube} ${\cal K}$ is a  subset of $\R^+ \times\R^N$, such that 
$\overline{\cal K}\cap ([0,t]\times\R^N)$ is a compact subset of $\R^{N+1}$ for any $ t\geq 0$.

\item
A set ${\cal K}\subset \R^+\times\R^N$ is {\it left lower semicontinuous} if
$$
\forall t>0,\; \forall x\in {\cal K}(t), \; {\rm if }\; t_n\to t^-, \;  {\rm then} \;
\exists  x_n\in{\cal K}(t_n) \;\mbox{\rm such that } x_n\to x\;.
$$

\item 
If $s= 1, 2$ or $(1,1)$,
a ${\cal{C}}^s$ {\it regular tube} ${\cal K}_r$ is a tube with a nonempty interior and whose 
boundary has a  ${\cal C}^s$ regularity, and is such that
at any point $(t,x)\in \partial{\cal K}_r,$ the outward unit normal 
$\nu_{(t,x)}^{{\cal K}_r}=(\nu_t,\nu_x)$ to ${\cal K}_r$ 
at $(t,x)$ satisfies 
\begin{eqnarray} \label{regular-tube}
\nu_x\neq0. 
\end{eqnarray}

\item  
The {\it normal velocity} $V^{{\cal K}_r}_{(t,x)}$ of a ${\cal C}^1$ regular tube ${\cal K}_r$ 
at the point $(t,x)\in\partial{\cal K}_r $ is defined by 
\begin{eqnarray} \label{def-vitesse}
V^{{\cal K}_r}_{(t,x)}=-\frac{\nu_t}{|\nu_x|},
\end{eqnarray}
where $\nu_{(t,x)}^{{\cal K}_r}= (\nu_t,\nu_x)$ is the outward unit normal to ${\cal K}_r$ at $(t,x)$.

\item
A ${\cal C}^1$ regular tube 
${\cal K}_r$ is {\it externally tangent} to a tube ${\cal K}$ at $(t,x)\in{\cal K}$ 
if
$$
{\cal K}\subset {\cal K}_r\;{\rm and}\; (t,x)\in\partial {\cal K}_r\;.
$$
It is {\it internally
tangent} to ${\cal K}$ at $(t,x)\in\widehat{{\cal K}}$ 
if $${\cal K}_r\subset {\cal K}\;{\rm and}\;(t,x)\in\partial {\cal K}_r\;.
$$

\item
We say that a sequence of ${\cal C}^{1,1}$ tubes $({\cal K}_n)$ converges to
some ${\cal C}^{1,1}$ tube ${\cal K}$ {\it in the ${\cal C}^{1,{\rm b}}$ sense}
if $({\cal K}_n)$ converges to ${\cal K}$ and $(\partial {\cal K}_n)$ converges to $\partial {\cal K}$
for the Hausdorff distance, and if there is an open neighborhood ${\cal O}$ of $\partial {\cal K}$
such that, if ${\bf d}_{\cal K}$ (respectively ${\bf d}_{{\cal K}_n}$) is the signed distance 
(\ref{def-dist-signee}) to ${\cal K}$ (respectively to ${\cal K}_n$), then $({\bf d}_{{\cal K}_n})$ and
$(D{\bf d}_{{\cal K}_n})$ converge uniformly to ${\bf d}_{\cal K}$ and $D{\bf d}_{\cal K}$
on ${\cal O}$ and $\| D^2{\bf d}_{{\cal K}_n} \|_\infty $ are uniformly bounded on ${\cal O}.$

\end{itemize}

\begin{Remark} \label{def-loc-t}  \rm $\;$  \\ 
1. The reason to introduce ${\cal C}^{1,1}$ and ${\cal C}^2$ tubes is clear when
looking at (\ref{fpp}) and (\ref{defh}): $\bar{h}$ is well defined if ${\cal K}_r$ is a ${\cal C}^{1,1}$
tube (see Section \ref{sec-propbarh}) and $F$ is well defined if ${\cal K}_r$ is a ${\cal C}^2$ tube 
(to be able to compute the curvature of ${\cal K}_r$). Therefore, (\ref{fpp}) is well defined when
${\cal K}_r$ is a ${\cal C}^2$ regular tube and, following the ideas of viscosity solutions (see \cite{cil92}),
${\cal C}^2$ regular tubes will play the role of ``test-functions" 
in the following definition. \\
2. For simplicity, we gave all the above definitions with tubes defined for all time $t\geq 0.$
But this is not the point, everything in the sequel is local in time, so we use the same 
definitions with tubes ${\cal K}$ where ${\cal K}(s)$ is only defined
in a neighborhood of some fixed $t.$ Note that smooth tubes will always be defined locally in time.
\end{Remark}

\begin{Definition} Let ${\cal K}$ be a tube and $K_0\in{\cal D}$ be an initial set.
\begin{enumerate}
\item
${\cal K}$ is {\rm a viscosity subsolution} to the
front propagation problem (in short FPP) (\ref{fpp}) 
if $\cal K$ is left lower semicontinuous and ${\cal K}(t)\in{\cal D}$ for any $t$,
and if, for any ${\cal C}^2$ regular
tube ${\cal K}_r$ externally tangent to ${\cal K}$ at some point $(t,x)$,
with ${\cal K}_r(t) \in{\cal D}$ and $t>0$, we have $$
V_{(t,x)}^{{\cal K}_r}\leq h_\lambda (x,{\cal K}_r(t))
$$
where $V_{(t,x)}^{{\cal K}_r}$ is the normal velocity of ${\cal K}_r$ at
$(t,x)$.

We say that ${\cal K}$ is a subsolution to the
FPP (\ref{fpp})  with initial position $K_0$ if 
${\cal K}$ is a subsolution and if $\overline{\cal K}(0)\subset \overline{K_0}$.

\item ${\cal K}$ is a {\rm viscosity supersolution} to the FPP (\ref{fpp}) 
if $\widehat{{\cal K}}$ is left lower semicontinuous, and ${\cal K}(t)\subset {\cal D}$
for any $t$, and if, for any ${\cal C}^2$ regular
tube ${\cal K}_r$ internally tangent to ${\cal K}$ at some point $(t,x)$, 
with ${\cal K}_r(t)\in{\cal D}$ and $t>0$, we have 
$$
V_{(t,x)}^{{\cal K}_r}\geq h_\lambda (x,{\cal K}_r(t))\;.
$$

We say that ${\cal K}$ is a supersolution to the
FPP (\ref{fpp}) with initial position $K_0$  if ${\cal K}$ is a supersolution
and if 
$\widehat{\cal K}(0)\subset \overline{\R^N\backslash K_0}$.

\item Finally, we say that a  tube ${\cal K}$ is a {\rm viscosity solution} to the front
propagation problem (with initial position $K_0$) if ${\cal K}$ is a sub- and a supersolution to the
FPP (with initial position $K_0$).
\end{enumerate}
\end{Definition}

\begin{Remark} \label{localisation-def} \rm 
The operator $h_\lambda$ defined in (\ref{defh}) is the sum of a local operator $F$ and a
nonlocal one $\bar{h}.$ As in the theory of viscosity solutions, 
we can localize arguments related to the local part of the operator. More precisely,
${\cal C}^2$ regularity of the boundary of the tube is required to compute
the curvature in $F,$ but only ${\cal C}^{1,1}$ regularity is needed to compute the
nonlocal part $\bar{h}.$ Therefore, the above definition is equivalent if we replace 
``for any ${\cal C}^2$ regular tube ${\cal K}_r$ internally (respectively externally) tangent 
to ${\cal K}$  at some point $(t,x)$...'' by 
``for any ${\cal C}^{1,1}$ regular tube ${\cal K}_r$ internally (respectively externally)
tangent to ${\cal K}$ at some point $(t,x)$ such that $\partial {\cal K}_r$ is ${\cal C}^2$
in a neighborhood of $(t,x)$...'' We will use this equivalent definition in the proof of Theorem
\ref{coconuts}.
\end{Remark}

\subsection{Regularity properties of the velocity $\bar{h}$}
\label{sec-propbarh}

We complete this part by recalling the regularity properties of the nonlocal 
term $\bar{h}$ defined by (\ref{Defh2}) and (\ref{pde}). 
These results were already given in \cite{cr04}, so we omit the proofs.
Here we assume that the set $S$ and the function $g$ satisfy assumptions (\ref{Hypfg}).

Because of the maximum principle, the function $\bar{h}$ is nonnegative and nondecreasing:
if $K_1\in{\cal D}$ and $K_2\in{\cal D}$ are closed and with a ${\cal C}^{1,1}$
boundary, if $K_1\subset K_2$ and if $x\in \partial K_1\cap K_2$,
then $0\leq \bar{h}(x,K_1)\leq \bar{h}(x,K_2).$

Furthermore, $\bar{h}$ is continuous in the following sense:
If $K_n$ and $K\in{\cal D}$  are closed subsets of $\R^N$ with ${\cal C}^{1,1}$
boundary such that $K_n$ converge to $K$ in the ${\cal{C}}^{1,{\rm b}}$ sense, 
if $x_n\in \partial K_n$ converge to 
$x\in \partial K$, then
$$
\lim_n \bar{h}(x_n,K_n)= \bar{h}(x,K).
$$
This is a straightforward application of \cite[Theorem 8.33]{gt83}. 

Next we give a result describing the behaviour of $\bar{h}$ for large ball:
\begin{Lemma} \label{growth}
For any $x_0\in\R^N$, there are constants $r_0>0$ and $\alpha>0$ such that
$$
\forall r\geq r_0,\; 
\forall x\in\partial B(x_0,r),\ \ \ 
\bar{h}(x,B(x_0,r))\leq \left\{\begin{array}{ll} 
\alpha r^{2-2N} & {\rm if }\; N\not= 2, \\
\displaystyle{\frac{\alpha}{r^2|\log(r)|^2}} & {\rm if }\; N=2.
\end{array}
\right. 
$$
Moreover, the constants $r_0$ and $\alpha$ only depend on 
$S$ and on $\|g\|_\infty$.
\end{Lemma}
\noindent
The proof is based on standard construction of supersolutions
to (\ref{pde}) for $\Omega= B(0,r)$, and so we omit it.

Lemma (\ref{growth}) states that $\bar{h}$ is small when $\Omega$ is a large ball.
On the contrary, the following lemma means that $\bar{h}$ is large
when ``$\Omega$ is close to $S.$'' 
For all $\gamma \geq 0,$ we introduce 
\begin{equation}\label{Sgamma}
S_\gamma =\{ x\in\R^N\;,\;  d_S(x)\leq \gamma\}.
\end{equation}
Then, we have
\begin{Lemma}\label{barhgrand}
There exist $\gamma_0>0$ and a constant $\alpha >0$ which depends only on $g$
and $S$
such that, for all $\gamma\in (0,\gamma_0),$
\begin{eqnarray*}
\bar{h}(x,S_\gamma )\geq \frac{\alpha}{\gamma^2}\quad \forall x\in \partial S_\gamma.
\end{eqnarray*}
\end{Lemma}

\noindent{\bf Proof of Lemma \ref{barhgrand}.}
Since $S$ has a ${\cal{C}}^2$ boundary, we can
fix $\gamma_0>0$ small enough such that ${\bf d}_S$ defined by (\ref{def-dist-signee})
is ${\cal{C}}^2$ in $S_{2\gamma_0}\backslash \{ {\bf d}_S < -2\gamma_0\}.$
We fix $\gamma\in (0,\gamma_0)$ and set $K = \{ {\bf d}_S \leq -\gamma \}.$ 
We note that $d_K={\bf d}_S+\gamma$ is ${\cal C}^2$
on $S_{2\gamma_0}\backslash \{ {\bf d}_S < -2\gamma_0\}.$
Moreover $d_{K}=\gamma$ on $\partial S$ and $d_{K}=2\gamma$ on $\partial S_\gamma.$
Set
\begin{equation}\label{plouic}
M=\max\{|\Delta d_K^2(x)|\; |\; 0\leq {\bf d}_S(x)\leq \gamma\}
\;{\rm and }\; m=\min\{g(x)\;|\; x\in\partial S\}\;.
\end{equation}
Finally we set $\Omega = S_\gamma$ and, for $\beta = {\rm e}^{-3M/4}m,$
we define
$$
\varphi (r)= \beta ({\rm e}^{-Mr/(4\gamma^2)} -1) \quad \forall r\in \R\;.
$$
We claim that
$$
u(x)=\varphi (d_K^2 (x)-(2\gamma)^2)
$$
is a subsolution of (\ref{pde}). Indeed, since $\varphi (0)=0,$
for all $x\in \partial S_\gamma,$ $u(x)=0.$ From the definition of $\beta,$
for all $x\in \partial S,$ $u(x)=\varphi (-3\gamma^2) \leq m\leq g.$
Setting $r_x =d_K^2 (x)-(2\gamma)^2,$ an easy computation gives
$$
-\Delta u(x) = -\varphi '' (r_x) |D(d_K^2)|^2 -\varphi '(r_x) 
\Delta (d_K^2).
$$
But $D(d_K^2)= 2 d_K Dd_K$ and $|Dd_K|=1.$ From (\ref{plouic}), we get
$$
-\Delta u(x) \leq -4\varphi ''(r_x)  d_K^2 +M |\varphi '(r_x)|.
$$  
A computation of the derivatives of $\varphi$ gives
$$
-\Delta u(x) \leq \frac{\beta M^2}{4\gamma^2} {\rm e}^{-Mr_x /(4\gamma^2)}
\left( 1-\frac{d_K^2 (x)}{\gamma^2}\right).
$$
For $x\in \Omega\backslash S,$ we have $d_K (x)\geq \gamma$ and therefore we obtain
$-\Delta u(x) \leq 0.$ Finally $u$ is a subsolution with $u\geq 0$ in
$\Omega$ and $u=0$ on $\partial S_\gamma.$ Thus, for all $x\in \partial S_\gamma,$
$$
\bar{h}(x, S_\gamma)\geq |Du (x)|^2 = \frac{M^2  {\rm e}^{-3M/2} m^2}{4\gamma^2}.
$$
\QED

We now recall the main regularity property of the map $\bar{h}$:

\begin{Lemma}\label{regul}
Let $R>0$ be some large constant and $\gamma>0$ be sufficiently small such that
$S_\gamma$ defined by (\ref{Sgamma}) 
has a ${\cal C}^2$ boundary. There is a constant $\theta>1/\gamma$ such that, for any 
compact set $K$  with ${\cal C}^{1,1}$ boundary such that $S_\gamma\subset {\rm int}(K)$ 
and $K\subset B(0,R-\gamma)$, for any $v\in \R^N$ with $|v|<1/\theta$ 
and any $x\in \partial K$, we have
\begin{equation}
\bar{h}(x+v,K+v)\geq (1-\theta|v|)^2\bar{h}(x,K).
\end{equation}
\end{Lemma}

\noindent For the proof, see \cite[ Proposition 2.4]{cr04}.


\section{Interposition theorems}

This part is devoted to interposition theorems in space and in space-time.
Such results are fondamental in the proof of the inclusion principle.
They play the same role as Jensen's maximum principle (see \cite{jensen88})
or Ishii's lemma (see \cite[Theorem 8.3]{cil92}) in the standard theory 
of viscosity solutions.

\subsection{An interposition theorem in $\R^N$}

Let us start with an interposition result for subsets of $\R^N$.
The following proposition is a direct consequence of Ilmanen interposition 
lemma \cite{ilmanen93} and can be found in \cite[Proposition 3.7]{cr04}. 

\begin{Proposition}[Interposition]\label{ilmanen} Let $K_1$ and $K_2$ be two closed subsets of 
$\R^N$, with $K_1$ compact and such that $K_1\subset\subset K_2$. Let $y_1\in K_1$ and 
$y_2\in \partial K_2$ be such that
$$
|y_1-y_2|=\min_{z_1\in K_1, z_2\in \partial K_2} |z_1-z_2|\;.
$$
Then there is some open subset $\Sigma_1$ of $\R^N$ with a ${\cal C}^{1,1}$ boundary, such that
$\Sigma_1$ is externally tangent to $K_1$ at $y_1$ (i.e., $K_1\subset \overline{\Sigma_1}$ and
$y_1\in\partial \Sigma_1$) and such that $\Sigma_2:=\Sigma_1+y_2-y_1$ is internally tangent
to $K_2$ at $y_2$ (i.e., $\Sigma_2\subset K_2$ and $y_2\in \partial \Sigma_2$).
\end{Proposition}
%
\begin{figure}[h]
\begin{center}
\epsfig{file=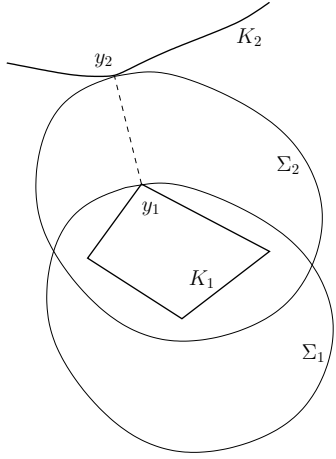, width=4.5cm} 
\end{center}
\caption{\label{dessin-ilmanen}
\textsl{Illustration of the result of Proposition \ref{ilmanen}.}}
\end{figure}
See Figure \ref{dessin-ilmanen} for an illustration of this proposition.
The key point in this result is that the smooth set $\Sigma_2$ 
internally tangent to $K_2$ is just a translation of the smooth set
$\Sigma_1$ externally tangent to $K_1$.

The ${\cal C}^{1,1}$ regularity of the sets $\Sigma_1$ and $\Sigma_2$ turns out to be optimal:
one cannot expect $\Sigma_1$ and $\Sigma_2$ to be ${\cal C}^2$ in general. 
Unfortunately the ${\cal C}^2$ regularity will be required in the sequel to be able to deal with 
curvature terms. 
In order to overcome this difficulty, one can approximate the sets $\Sigma_1$ and $\Sigma_2$ 
in the following way:

\begin{Theorem}[Approximation]\label{app} Let $K_1$, $K_2$, $y_1$, $y_2$, $\Sigma_1$ and $\Sigma_2$ be as
in Proposition \ref{ilmanen} and $\delta>0$ be sufficiently small. Then 
there exists $\Sigma_{1,n}$ and $\Sigma_{2,n}$ open subsets of $\R^N$
with ${\cal C}^{1,1}$ boundary, converging
respectively to $\Sigma_1$ and $\Sigma_2$ in the ${\cal C}^{1,{\rm b}}$ sense, there exists $y_{1,n}\in K_1$
and $y_{2,n}\in \partial K_2$ converging respectively to $y_1$ and $y_2$, 
and there exists $(N-1)\times (N-1)$ matrices $X_1$, $X_2$ such that
\begin{itemize}
\item[(i)] $\Sigma_{1,n}$ is externally tangent to $K_1$ at $y_{1,n}$
and $\Sigma_{2,n}$ is internally tangent to $K_2$ at $y_{2,n}$.

\item[(ii)] For $i=1$ and $2,$ 
$\Sigma_{i,n}$ is of class ${\cal C}^2$ in a neighbourhood of $y_{i,n}$
with $\displaystyle{ \lim_n H^{\Sigma_{i,n}}_{y_{i,n}}\to X_i,  }$ and
\begin{equation}\label{AAaa}
-\frac{1}{\delta} I_{2(N-1)}\leq 
\left(\begin{array}{cc}
X_1 & 0\\ 0 & -X_2
\end{array}\right)\leq
\frac{1}\delta \left(\begin{array}{cc}
I_{N-1} & -I_{N-1}\\ -I_{N-1} & I_{N-1}
\end{array}\right)\;.
\end{equation}
\end{itemize}
\end{Theorem}

\begin{Remark} \rm $\;$  \\
1. Note carefully that the two approximations are not independent because of the inequalities (\ref{AAaa}),
which implies in particular that $X_1\leq X_2$.\\
2. By ``$\delta>0$ sufficiently small", we mean $\delta\in(0, |y_1-y_2|/(2+|y_1-y_2|))$.
\end{Remark}

\noindent
The proof of Theorem \ref{app} is very similar to the (more difficult) proof of the second part 
of Theorem \ref{interposition} below, so we omit it. 

\subsection{Interposition by regular tubes}\label{InterpoTube}

The aim of this part is to extend previous results for subsets $\Sigma_1, \Sigma_2\subset \R\times\R^N$ 
which are regular tubes. The point here is to be able to construct tubes satisfying the
regularity assumption (\ref{regular-tube}). 

For this we introduce some notations. 
In $\R\times \R^N$ we work with the norm (where $\sigma>0$ is fixed)
$$
|(t,x)|_\sigma=\left(\frac{1}{\sigma^2}t^2+|x|^2\right)^\frac12\;.
$$
For any subset $E$ of $\R^{N+1}$, we note the distance to $E$ for this norm
$$
d^\sigma_E(t,x)=\inf_{(s,y)\in E} |(s,y)-(t,x)|_\sigma\;.
$$
For any two subsets $A_1$ and $A_2$ of $\R^{N+1}$, we define the minimal
distance between $A_1$ and $A_2$ by
$$
e(A_1,A_2)=\inf_{(t_1,x_1)\in A_1, \; (t_2,x_2)\in A_2} |(t_2,x_2)-(t_1,x_1)|_\sigma\;.
$$

We consider the following transversality condition:
\begin{eqnarray} \label{Transversality}
\begin{array}{c}
{\rm for} \ C_1\subset\subset C_2 \subset \R\times \R^N \ 
{\rm with} \ C_1 \ {\rm compact \ and} \ C_2 \ {\rm closed}, \\ 
{\rm and \ for \ any} \  (\bar{s}_1,\bar{y}_1)\in C_1 \ {\rm and \ any} \ (\bar{s}_2,\bar{y}_2)\in 
\widehat{C_2}, \\
{\rm if }\; |(\bar{s}_1,\bar{y}_1)-(\bar{s}_2,\bar{y}_2)|_\sigma=e(C_1,\widehat{C_2}), \ {\rm then} \
\bar{s}_1>0, \bar{s}_2>0\;{\rm and }\;
\bar{y}_1\neq \bar{y}_2.
\end{array}
\end{eqnarray}

\begin{Theorem}\label{interposition} Let $C_1$
and $C_2$ be such that (\ref{Transversality}) holds. 
Let us fix $(\bar{s}_1,\bar{y}_1)\in C_1$ and $(\bar{s}_2,\bar{y}_2)\in \widehat{C_2}$ with
$$
|(\bar{s}_1,\bar{y}_1)-(\bar{s}_2,\bar{y}_2)|_\sigma=e(C_1,\widehat{C_2})\;.
$$
\begin{enumerate}
\item {\bf Interposition:}
There exists a ${\cal C}^{1,1}$ regular tube $\Sigma_1$, defined on an open interval $I$ (see Remark \ref{def-loc-t}.2),
such that $\Sigma_1$ is externally tangent to $C_1$ at $(\bar{s}_1,\bar{y}_1)$, with $\bar{s}_1\in I$, 
and $\Sigma_2:=\Sigma_1 +(\bar{s}_2,\bar{y}_2)-(\bar{s}_1,\bar{y}_1)$ is internally tangent to $C_2$ at $(\bar{s}_2,\bar{y}_2)$.

\item {\bf Joint approximation by ${\cal C}^2$ tubes:}
Futhermore, for any $\delta>0$ sufficiently small, there exists ${\cal C}^{1,1}$ regular tubes $\Sigma_{1,n}$ 
and $\Sigma_{2,n}$ converging
respectively to $\Sigma_1$ and $\Sigma_2$  in the ${\cal C}^{1,{\rm b}}$ sense, 
there exists $(\bar{s}_{1,n}, \bar{y}_{1,n})\in C_1$
and $(\bar{s}_{2,n}, \bar{y}_{2,n})\in \widehat{C_2}$ converging respectively to $(\bar{s}_1,\bar{y}_1)$ and $(\bar{s}_2,\bar{y}_2)$, 
and there exists $(N-1)\times (N-1)$ matrices $X_1$, $X_2$ such that
\begin{itemize}
\item[(i)] $\Sigma_{1,n}$ is externally tangent to $C_1$ at $(\bar{s}_{1,n}, \bar{y}_{1,n})$
and $\Sigma_{2,n}$ is internally tangent to $C_2$ at $(\bar{s}_{2,n},\bar{y}_{2,n})$.

\item[(ii)] For $i=1$ and $2,$ $\Sigma_{i,n}$ is of class ${\cal C}^2$ in a neighbourhood of $(\bar{s}_{i,n},\bar{y}_{i,n})$
with $\displaystyle{ \lim_n H^{\Sigma_{i,n}(\bar{s}_{i,n})}_{\bar{y}_{i,n}}\to X_i}$ and
\begin{equation}\label{AA}
-\frac{1}{\delta} I_{2(N-1)}\leq 
\left(\begin{array}{cc}
X_1 & 0\\ 0 & -X_2
\end{array}\right)\leq
\frac{1}\delta \left(\begin{array}{cc}
I_{N-1} & -I_{N-1}\\ -I_{N-1} & I_{N-1}
\end{array}\right)\;.
\end{equation}
\end{itemize}
\end{enumerate}
\end{Theorem}
The proof of this theorem is done in Section \ref{preuve-thm-interposition}.

\begin{Remark} \rm $\;$  \\
1. Inequality (\ref{AA}) implies that $X_1\leq X_2$. Although we only use this latter inequality in the sequel, inequality
(\ref{AA}) allows to treat equations with $F$ depending on $x$ (see for instance \cite{cil92}). Let us once 
again point out that the two approximations are not independent because of (\ref{AA}).\\
2. Thanks to the ${\cal C}^{1,{\rm b}}$ convergence of $\Sigma_{1,n}$ and $\Sigma_{2,n}$
to $\Sigma_1$ and $\Sigma_2$ respectively, one also has:
\begin{equation}\label{nunu}
\lim_n\nu^{\Sigma_{1,n}(\bar{s}_{1,n})}_{\bar{y}_{1,n}}=\lim_n \nu^{\Sigma_{2,n}(\bar{s}_{2,n})}_{\bar{y}_{2,n}}= 
\nu^{\Sigma_1(\bar{s}_1)}_{\bar{y}_1}=\nu^{\Sigma_2(\bar{s}_2)}_{\bar{y}_2}, 
\end{equation} 
\begin{equation}\label{VV}
\lim_n V^{\Sigma_{1,n}}_{(\bar{s}_{1,n},\bar{y}_{1,n})}=\lim_n V^{\Sigma_{2,n}}_{(\bar{s}_{2,n},\bar{y}_{2,n})}= 
V^{\Sigma_1}_{(\bar{s}_1,\bar{y}_1)}=V^{\Sigma_2}_{(\bar{s}_2,\bar{y}_2)}\;.
\end{equation} 
3. By ``$\delta>0$ sufficiently small, we mean: $\delta \in (0, e(C_1,\widehat{C_2})/(2+e(C_1,\widehat{C_2})))$.
\end{Remark}

\subsection{Existence of the regular interposition tubes}

Let us introduce a new notation: 
if $\Sigma$ is a tube defined on some open interval $I$ (see Remark \ref{def-loc-t}.2), 
then we set
\begin{equation}\label{DefBd}
{\rm bd}(\Sigma):=\bigcup_{t\in I} \partial \overline{\Sigma}(t) \;.
\end{equation}
The following result is the key point in the proof of the existence of the regular
interposition tubes of Theorem \ref{interposition} part (i).

\begin{Proposition}\label{troistrois} Let $C_1$, $C_2$, $(\bar{s}_1,\bar{y}_1)\in C_1$, $(\bar{s}_2,\bar{y}_2)\in C_2$ be
as in Theorem \ref{interposition}.
There exist a ${\cal C}^{1,1}$ regular tube $\Sigma$ defined on some interval $I$ and some 
$(t,x)\in](\bar{s}_1,\bar{y}_1),(\bar{s}_2,\bar{y}_2)[$ such that
\begin{equation}\label{EqDesiree}
t\in I, \quad
x\in \partial \Sigma(t)\quad {\rm and }\quad 
e (C_1\;,\; \widehat{C_2})=
e (C_1\;,\; {\rm bd}( \Sigma))+
e (\Sigma\;,\; \widehat{C_2})\;.
\end{equation}
\end{Proposition}
Above, $](\bar{s}_1,\bar{y}_1),(\bar{s}_2,\bar{y}_2)[ $ denotes the open segment joining 
$(\bar{s}_1,\bar{y}_1)$ and $(\bar{s}_2,\bar{y}_2)$.\\

\noindent{\bf Proof of Proposition \ref{troistrois}.} 
Let us first fix some notation needed throughout the proof: 
we set 
\begin{itemize}
\item $\bar{e} := e( C_1\;,\; \widehat{C_2})$,
\item $\displaystyle{   E := \{ (t,x)\in \R^+\times \R^N  :  (t,x)\in ] (s_1,y_1) , (s_2,y_2) [ \ {\rm where} }$\\ 
$\displaystyle{  \hspace{1cm} (s_1,y_1)\in C_1 \ {\rm and} \ (s_2,y_2) \in \widehat{C_2} \ {\rm satisfy} 
\ |(s_1,y_1)-(s_2,y_2)|_\sigma=\bar{e}  \}, }$
\item $\displaystyle{   
A_\rho:=\{(t,x)\in \R^+\times \R^N\; :\; d^\sigma_{C_1} (t,x)>\rho\;{\rm and }\; d^\sigma_{\widehat{C_2}}(t,x)>\rho\}
\ {\rm for}\; \rho\in(0,\bar{e}/2)\;, }$
and, if $I$ is an interval, then $A_\rho(I)=A_\rho\cap (I\times \R^N)$,
\item $\displaystyle{ (t(s),x(s)):=s(\bar{s}_1, \bar{y}_1)+(1-s)(\bar{s}_2, \bar{y}_2)}$ for all $s\in(0,1)$
and 
\begin{eqnarray}\label{doubleetoile}
(\bar{t},\bar{x}):=(t(1/2),x(1/2)),
\end{eqnarray}
\item $I_\tau:=(\bar{t}-\tau,\bar{t}+\tau)$ for all $\tau>0.$
\end{itemize}
For later use we note that
\begin{equation}\label{tsxs}
d^\sigma_{\widehat{C_2}}(t(s),x(s))= s\bar{e}\quad{\rm and }\quad d^\sigma_{C_1}(t(s),x(s))=(1-s)\bar{e}\quad
{\rm for \ all} \ s\in(0,1)\;,
\end{equation}
because of the definition of $(\bar{s}_1, \bar{y}_1)$ and $(\bar{s}_2, \bar{y}_2)$.
Moreover, for a point $(t,x)\in \R^+\times \R^N$,
the equality $\bar{e}- d^\sigma_{\widehat{C_2}} (t,x)=d^\sigma_{C_1} (t,x)$ holds if and only if $(t,x)\in \overline{E}$.

We now reduce the construction of the tube $\Sigma$ to the construction of a 
suitable function $w$:

\begin{Lemma}\label{ineq} 
Let $I$ be a nonempty open interval of $\R^+$ and $w: A_\rho(I)\to\R$ be 
of class ${\cal C}^{1,1}$ (for some $\rho\in (0,\bar{e}/2)$) and such that 
\begin{equation}\label{blabla}
\bar{e}- d^\sigma_{\widehat{C_2}} (s,y) \leq w (s,y) \leq d^\sigma_{C_1} (s,y)\qquad \forall (s,y)\in A_\rho(I)\;.
\end{equation}
We also assume that there is some $\gamma\in(\rho,\bar{e}-\rho)$ and some $(t,x)\in E\cap A_\rho(I)$ with $w(t,x)=\gamma$ and such that
\begin{equation}\label{GradNonZero}
D_x w(s,y)\neq 0\qquad \forall (s,y)\in A_\rho(I)\quad{\rm with}\quad w(s,y)=\gamma\;.
\end{equation}
Then the set $\Sigma=\{(s,y)\in A_\rho(I)\;|\; w(s,y)\leq \gamma\}$
satisfies the requirements of Proposition \ref{troistrois}.
\end{Lemma}

\noindent{\bf Proof of Lemma  \ref{ineq}.} Let us first check that $\Sigma$ is a tube of class ${\cal C}^{1,1}$ in 
the intervall $I$. Because of assumption (\ref{GradNonZero}) it suffices to show that
${\rm bd}(\Sigma)\subset A_\rho(I)$ where ${\rm bd}(\Sigma)$ is defined by (\ref{DefBd}).
Using (\ref{blabla}) and the fact that  $\gamma\in(\rho,\bar{e}-\rho)$, we have,
for all $(s,y)\in {\rm bd}(\Sigma)$,
$$
d^\sigma_{\widehat{C_2}} (s,y)\geq \bar{e} -\gamma >\rho\quad{\rm and }\quad
d^\sigma_{C_1} (s,y)\geq \gamma>\rho\;.
$$
Hence $(s,y)\in A_\rho(I)$.

Finally we show that (\ref{EqDesiree}) holds. If $w(s,y)=\gamma$, then 
$d^\sigma_{C_1}(s,y)\geq \gamma$ while $d^\sigma_{\widehat{C_2}}(s,y)\geq \bar{e}-\gamma$.
Hence $e(C_1, {\rm bd}(\Sigma))\geq \gamma,$ $e(\widehat{C_2}, \Sigma)\geq \bar{e}-\gamma$ and
$$
e(C_1, {\rm bd}(\Sigma))+e(\widehat{C_2}, \Sigma)\geq e(C_1,\widehat{C_2})\;.
$$
For the reverse inequality let us first recall that $\bar{e}- d^\sigma_{\widehat{C_2}} (t,x)=d^\sigma_{C_1} (t,x)$
because $(t,x)\in E$. This implies that $\gamma=d^\sigma_{C_1}(t,x)\geq e(C_1, {\rm bd}(\Sigma))$ and 
$\bar{e}-\gamma= d^\sigma_{\widehat{C_2}}(t,x)\geq e(\widehat{C_2}, \Sigma)$.
Hence (\ref{EqDesiree}) holds. 
\QED

\begin{figure}[h]
\begin{center}
\epsfig{file=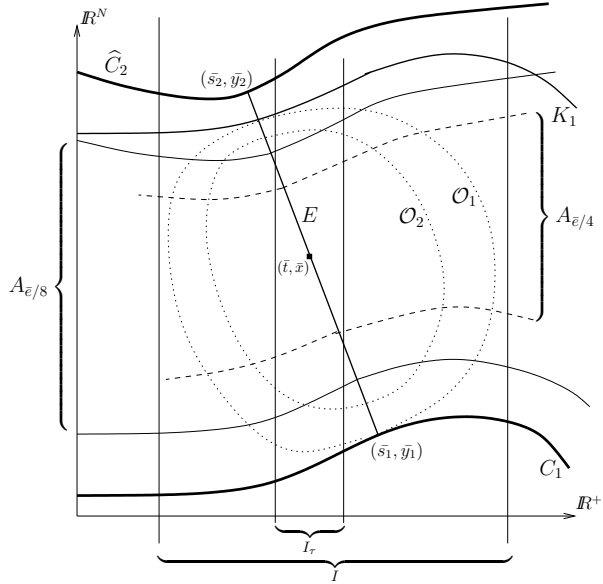, width=8cm} 
\end{center}
\vspace*{-0.7cm}
\caption{\label{dessin-interregu}
\textsl{Illustration of the proof of Proposition \ref{troistrois}.}}
\end{figure}

Next we turn to the construction of a function $w$ satisfying the assumptions of
Lemma \ref{ineq}. We advice the reader to look at Figure \ref{dessin-interregu} to follow the rest
of the proof of Proposition \ref{troistrois}.

The first step is the following result given in \cite{cardaliaguet01}: let us set 
$$
K_1 :=\{ (t,x)\in \R^+\times \R^N : d^\sigma_{C_1} (t,x) \leq 15\bar{e}/16 \}\;.
$$
Then, we have
\begin{Lemma}
\cite[Lemma 5.1]{cardaliaguet01}\label{carda}
The function $\varphi (t,x)= d_{\partial K_1}^\sigma (t,x)$ is ${\cal{C}}^{1,1}$ in a bounded
open neighborhood
${\cal{O}}_1$ of $E\cap (K_1 \backslash \partial K_1)$ and 
\begin{eqnarray}\label{coincement}
\bar{e}- d^\sigma_{\widehat{C_2}} (t,x) \leq \frac{15 \bar{e}}{16} -\varphi (t,x) \leq d^\sigma_{C_1} (t,x)
\quad \forall (t,x)\in \R^+\times\R^N\;.
\end{eqnarray}
\end{Lemma}

We now show that $\varphi$ has a nonvanishing spatial gradient in a neighborhood
of $E\cap {\cal O}_1$.
Let $(t,x)\in E\cap {\cal{O}}_1.$ Then there exists $(s_1,y_1)\in C_1$, $(s_2,y_2) \in \widehat{C_2}$ 
such that $|(s_1,y_1)-(s_2,y_2)|_\sigma=\bar{e}$ and $(t,x)\in ] (s_1,y_1) , (s_2,y_2) [.$
Therefore, $D \varphi (t,x) = ((s_2,y_2)-(s_1,y_1)))/\bar{e}$. 
From assumption (\ref{Transversality}), we know that $y_1\not= y_2.$ Thus
$D_x \varphi (t,x) \not=0$ for all  $(t,x)\in E\cap {\cal{O}}_1.$
By continuity of $D_x\varphi$ in ${\cal{O}}_1,$
there is an open set ${\cal O}_2\subset {\cal O}_1$ which contains $E\cap A_{\bar{e}/8}$
and such that 
\begin{eqnarray} \label{nonzero1}
\eta := \mathop{\rm min}_{(t,x)\in \overline{{\cal{O}}_2}} | D_x\varphi (t,x) | > 0.
\end{eqnarray}

We are now going to modify $\varphi$ far away from $E$. 
For this we need a technical lemma:

\begin{Lemma} \label{lemme-tech2}
Let $f,g,h$ be continuous functions in $\R^k$
such that $g\leq f\leq h$ in $\R^k.$ Suppose that $K:=\{h=g\}$ is non empty and compact and
that there is some open neighbourhood $U$ of $K$ such that 
$f$ is ${\cal C}^{1,1}$ in $U$. 

Then, for any $\tilde{\eta} >0$ and for any open subset $U'$ such that
$K\subset U'\subset\subset U$, there is a function $\psi : \R^k\to \R$ such that: 
\begin{itemize}
\item[(i)] $\psi$ is ${\cal{C}}^{1,1}_{loc}$ in $\R^k$ and $\psi$ is ${\cal{C}}^\infty$ in $\R^k\backslash U',$ 
\item[(ii)] $g\leq \psi\leq h$ in $\R^k$ and $g<\psi<h$ in $\R^k\backslash U',$
\item[(iii)] $|D(\psi-f)| \leq\tilde{\eta}$ in $\overline{U'}.$
\end{itemize}
\end{Lemma}

\noindent{\bf Proof of Lemma \ref{lemme-tech2}.} Let $U_1$ and $U_2$ be two open subsets of $\R^k$
such that $K\subset U_2\subset\subset U_1\subset\subset U'$ and fix some smooth map
$\theta:\R^k\to [0,1]$ such that $\theta=1$ in $\overline{U_2}$ and $\theta=0$ in $\R^k\backslash U_1.$
Then we consider a smooth map $\xi:\R^k\to\R$ such that 
$g<\xi<h$ in $\R^k\backslash U_2$,
$$
\|\xi-f\|_{L^\infty(U'\backslash U_2)}\leq \frac{\tilde{\eta}}{2\|D\theta\|_\infty} \quad{\rm and }\quad
\|D\xi-Df\|_{L^\infty(U'\backslash U_2)}\leq \frac{\tilde{\eta}}{2}\;.
$$
The construction of such a function $\xi$ is possible because $g<f<h$
outside of $K$ and $f$ is ${\cal C}^{1,1}$ in $U$.
Then we set 
$$
\psi(x)=\theta(x)f(x)+(1-\theta(x)) \xi(x)\qquad \forall x\in\R^k\;.
$$
Note first that (i) and (ii) obviously hold. As for (iii), it clearly holds in $U_2$ since $\psi=f$
in $U_2$. Moreover, for  $x\in U'\backslash U_2$, we have
$$
|D(\psi-f)(x)|\leq |D(\xi-f)(x)|+|f(x)-\xi(x)||D\theta(x)|\leq \frac{\tilde{\eta}}{2}+
\frac{\tilde{\eta} |D\theta(x)|}{2\|D\theta\|_\infty}\leq \tilde{\eta}\;.
$$
\QED

Next, we apply Lemma \ref{lemme-tech2} with
$k:=N+1$, $U:={\cal O}_1$, $U':={\cal O}_2$, $\tilde{\eta}:=\eta/2$ where $\eta$
is given by (\ref{nonzero1}),
$$
g(t,x):=\bar{e}-d^\sigma_{\widehat{C_2}}(t,x)-d_{A_{\bar{e}/8}}^\sigma (t,x),\qquad 
h(t,x):=d^\sigma_{C_1}(t,x)+d_{A_{\bar{e}/8}}^\sigma (t,x)
$$
and $f(t,x):=15\bar{e}/16-\varphi(t,x)$ for $(t,x)\in\R^+\times\R^N$. We extend $f$, $g$ and $h$ for $t\leq0$ 
by setting $f(t,x)=f(0,x)$, $g(t,x)=g(0,x)$ and $h(t,x)=h(0,x)$.
From Lemma \ref{carda} we have $g\leq f\leq h$ in $\R^{N+1}$.
Moreover, from assumption (\ref{Transversality}), the set $K:= \{g=h\}=E\cap \overline{A_{\bar{e}/8}}$
is compact and contained in $(0,+\infty)\times \R^N$
and we know from Lemma \ref{carda} that $f$ is ${\cal C}^{1,1}$ in the neighbourhood ${\cal O}_1$ of $K$. 
Lemma \ref{lemme-tech2} states that there is a map $\psi: \R^{N+1}\to\R$ such that
\begin{equation}\label{PptePsi}
\left\{\begin{array}{l}
\mbox{\rm (i) $\psi$ is ${\cal{C}}^{1,1}_{loc}$ in $\R^{N+1}$ and ${\cal{C}}^\infty$ in $\R^{N+1}\backslash {\cal O}_2,$ }\\
\mbox{\rm (ii) $\bar{e}-d^\sigma_{\widehat{C_2}}\leq \psi\leq d^\sigma_{C_1}$ in $A_{\bar{e}/8},$}\\
\mbox{\rm (iii) $\bar{e}-d^\sigma_{\widehat{C_2}}<\psi<d^\sigma_{C_1}$ in $A_{\bar{e}/8}\backslash {\cal O}_2,$}\\
\mbox{\rm (iv) $|D(\psi-f)| \leq \eta/2$ in $\overline{{\cal O}_2}.$}
\end{array}\right.
\end{equation}
Putting together (\ref{nonzero1}) and (iv) implies that
$$
|D_x\psi(t,x)|\geq |D_xf(t,x)|-|D_x(\psi-f)(t,x)|\geq \frac{\eta}{2}\quad \forall (t,x)\in \overline{{\cal O}_2}\;.
$$
We now choose two open subsets $U_1$ and $U_2$ of $\R^N$ and some $\tau>0$
such that ${\cal O}_2(\bar{t})\subset\subset U_2\subset\subset U_1$ (recall that $\bar{t}$ is
defined by (\ref{doubleetoile})) and
\begin{equation}\label{IneqDxPsi}
|D_x\psi(t,x)|\geq \frac{\eta}{4}\qquad \forall (t,x)\in I_\tau\times \overline{U_1}\;.
\end{equation}
We also fix some smooth function $\theta:\R^N\to[0,1]$ such that $\theta=1$
in $U_2$, $\theta=0$ in $\R^N\backslash U_1$ and we set
$$
w(t,x)=\theta(x)\psi(t,x)+(1-\theta(x)) \psi(\bar{t},x)\qquad \forall (t,x)\in\R^{N+1}\;.
$$
Note that $w$ belongs to ${\cal C}^{1,1}_{loc}(\R^{N+1})\cap {\cal C}^\infty(\R\times (\R^N\backslash U_1))$.
We claim that we can choose
$\tau>0$ sufficiently small such that 
\begin{equation}\label{ChoixTau}
\left\{\begin{array}{l}
\mbox{\rm (i) $\bar{e}-d^\sigma_{\widehat{C_2}}\leq w\leq d^\sigma_{C_1}$ in $A_{\bar{e}/4}(I_\tau),$}\\
\mbox{\rm  (ii) $|D_x w|\geq \eta/8$ in $I_\tau\times \overline{U_1}.$}
\end{array}\right.
\end{equation}
Let us prove the first assertion. Let $(t,x)\in A_{\bar{e}/4}(I_\tau).$ On the one hand, if 
$(t,x)\in A_{\bar{e}/4}(I_\tau)\cap (I_\tau \times U_2),$ then $\theta (x)=1$ and $w(t,x)=\psi (t,x).$
Since $A_{\bar{e}/4}(I_\tau)\subset  A_{\bar{e}/8},$ we conclude from (\ref{PptePsi}(ii)).
On the other hand, suppose $(t,x)\in A_{\bar{e}/4}(I_\tau)\backslash (I_\tau \times U_2).$
In particular, $(t,x)\notin {\cal O}_2(\bar{t})$ and $(\bar{t},x)\notin {\cal O}_2.$ Therefore, from
(\ref{PptePsi}(iii)), we have
$$
\bar{e}-d^\sigma_{\widehat{C_2}}(\bar{t},x)\;<\;\psi(\bar{t},x) \;<\;d^\sigma_{C_1}(\bar{t},x).
$$
Using the uniform continuity of the above distance functions in the compact set 
$\overline{A_{\bar{e}/4}(I_\tau)},$ we obtain that, for $\tau$ small enough, 
$$
\bar{e}-d^\sigma_{\widehat{C_2}}(t,x)\;<\;\psi(\bar{t},x) \;<\;d^\sigma_{C_1}(t,x) \quad
\forall (t,x)\in A_{\bar{e}/4}(I_\tau)\backslash (I_\tau\times U_2)\;.
$$
Combining with (\ref{PptePsi}(ii)), we conclude also in this case.

For the second assertion,
we notice that, for $(t,x)\in I_\tau\times \overline{U_1}$, we have
$$
|D_x (w-\psi)(t,x)|\leq |1-\theta (x)| |D_x(\psi(t,x)-\psi(\bar{t},x))|
+|\psi(t,x)-\psi(\bar{t},x)||D\theta(x)|
$$
with a right-handside smaller than $\eta/8$ provided $\tau$ is sufficiently small, 
because $\psi\in{\cal C}^{1,1}_{loc}$. 
Then, for
$(t,x)\in I_\tau\times \overline{U_1}$, we have, from the choice of $U_1$ and $\tau$ in (\ref{IneqDxPsi}),
$$
|D_x w(t,x)|\geq |D_x\psi(t,x)|-|D_x (w-\psi)(t,x)|\geq \eta/8\;,
$$
which proves the second statement. \\

We now fix $\sigma\in(0,\bar{e}/4)$ such that 
$$
(t(s),x(s))\in I_\tau\times {\cal O}_2(\bar{t}) \qquad \forall s\in (1/2-\sigma, 1/2+\sigma)\;.
$$
This is possible because $(\bar{t},\bar{x})=(t(1/2),x(1/2))$ belongs to ${\cal O}_2(\bar{t})$.
From (\ref{tsxs}), an easy calculation gives $w(t(s),x(s))=(1-s)\bar{e}.$
Since $w$ is smooth in $\R\times (\R^N\backslash \overline{U_1})$, Sard Lemma states that we can find
a level $\gamma\in ((1/2-\sigma)\bar{e}, (1/2+\sigma)\bar{e})$ such that $\gamma$
is a non critical value of $w$ in $\R\times (\R^N\backslash U_1)$. 
We claim that $w$ and $\gamma$ satisfy the requirements of Lemma \ref{ineq}. Note first that, 
for $s=(\bar{e}-\gamma)/\bar{e}$,
 the point $(t(s), x(s))$ belongs to $E\cap A_{\bar{e}/4}(I_\tau)$ and satisfies $w(t(s), x(s)) =\gamma$.
Moreover, (\ref{blabla}) holds from (\ref{ChoixTau}(i)). Finally we show that 
(\ref{GradNonZero}) holds. Indeed, if $w(t,x)=\gamma$ for some $(t,x)\in  A_{\bar{e}/4}(I_\tau)$, then 
either $x\in \overline{U_1}$, in which case $D_xw(t,x)\neq 0$ thanks to (\ref{ChoixTau}(ii)), 
or $x\notin \overline{U_1}$
and $Dw(t,x)=(0,D_xw(,tx))\neq0$ because $\gamma$ is a non critical value of $w$ 
in $I_\tau\times (\R^N\backslash \overline{U_1})$. 
In each case we have $D_xw(t,x)\neq0$. This completes the proof of Proposition \ref{troistrois}.
\QED

\subsection{${\cal C}^2$ regularization of a ${\cal C}^{1,1}$ tangent surface near contact
points}

The aim of this section is to show the following fact:
if a ${\cal{C}}^{1,1}$ surface $\Sigma$ is externally tangent to a set $K$ at a point
$y,$ then it is possible to find a ${\cal C}^{1,1}$ surface $\tilde{\Sigma}$
which is close to $\Sigma$ (in the ${\cal C}^{1,{\rm b}}$ sense) and is still externally tangent
to $K$ at a point $\tilde{y}$ close to $y.$ Moreover, $\tilde{\Sigma}$ is more regular than
$\Sigma,$ namely is ${\cal{C}}^2$ in a neighborhood of $\tilde{y}.$ In particular, we can use
$\tilde{\Sigma}$ as a test set to estimate the curvature (see Remark \ref{localisation-def}). 

Let us give the exact assumptions:
\begin{eqnarray}\mbox{
\begin{minipage}{11.5cm} \label{assumA}
\noindent Let $K$ be a subset of $\R^{{k}}$ for ${{k}}\geq 1$
and $\Sigma$ be an open set with a ${\cal C}^{1,1}$ boundary $\partial \Sigma$,
which is externally tangent to $K$ at some point $y\in\partial K$. 
Let $x\notin K$ be such that $y$ is the unique projection of $x$ onto $K$ and $p:=D d_K(x)$
is the outward normal to $\Sigma$ at $y$. Suppose
that there is a sequence of points $x_n\to x$, where $d_K$ is twice differentiable
with first and second derivative denoted respectively $p_n$ and $X_n$, 
and finally assume that $p_n$ converges to $p$ while $X_n$ converges to some $X$. 
\end{minipage}
}
\end{eqnarray}
\noindent
Note that, by usual properties of the distance function at differentiability points, 
the projection of $x_n$ onto $K$ is unique and converges to $y$.
We denote by $y_n$ this projection.

\begin{Proposition} \label{turlututu}
Under Assumption (\ref{assumA})
we can find a sequence of open sets $\Sigma_n$ with ${\cal C}^{1,1}$ boundary such that
\begin{itemize}
\item[(i)] $\Sigma_n$ is externally tangent to $K$ at $y_n$,

\item[(ii)] $\Sigma_n$ has a ${\cal C}^2$ boundary in a neighbourhood of $y_n$, with normal $p_n$
and curvature equal to the restriction to $(p_n)^\perp$ of $-(X_n -\frac{1}{n} I_{{k}})$,

\item[(iii)] $\Sigma_n$ converges to $\Sigma$ in the ${\cal C}^{1,{\rm b}}$ sense. 

\end{itemize}
\end{Proposition}
Before starting the proof of the proposition, we need two lemmas. The first one builds, from
the derivatives of the distance function at a point $a,$ a map $\phi$ which has a local maximum on $K$
at the point $b$, projection of $a$ onto $K$:

\begin{Lemma}\label{pouet1} Suppose that $a\notin K$ and that $d_K$ is twice differentiable at $a$. Let 
$b$ be the projection of $a$ onto $K$. Then, for any $\alpha>0$, the (smooth) function
$$
\phi(z)= \langle Dd_K(a),z-b\rangle+\frac12 \langle (D^2d_K(a)-\alpha I_{{k}})(z-b), z-b\rangle
$$
has a strict local maximum at $b$ on $K$.
\end{Lemma}

\noindent {\bf Proof of Lemma \ref{pouet1}.} For any $z\in K$, we have
$$
\phi(z)= d_K(z+a-b)-\frac{\alpha}{2} |z-b|^2-d_K(a)-|z-b|^2\epsilon(z-b)
$$
because $d_K$ is has a second order Taylor expansion at $a$. But, since $z\in K$, we have
$d_K(z+a-b) \leq |(z+a-b )-z|= |a-b|=d_K(a)$. Therefore
$$
\phi(z)\leq -\frac{\alpha}{2} |z-b|^2-|z-b|^2\epsilon(z-b)
$$
which is negative as soon as $z\in K$ is sufficiently close to $b$ and $z\neq b$.
\QED

From now on, we fix a smooth function $\theta:\R^+\to [0,1]$ such that
$\theta$ is nonincreasing, $\theta=1$ on $[0,1/2]$ and
$\theta=0$ on $[1,+\infty)$.

We will use several times below the following interpolation. The proof relies on
straightforward computations, so we skip it.

\begin{Lemma}\label{pouet2} Let $\phi$ and $\psi$ some ${\cal C}^{1,1}$ functions in some open set ${\cal O}$.
Let $\bar{y}\in{\cal O}$ be such that $\phi(\bar{y})=\psi(\bar{y})$, and let us set, for any $\rho>0$, 
$$
\xi_\rho(z)= \phi(z)\theta_\rho(z)+\psi(z)(1-\theta_\rho(z))\qquad 
{where } \ \theta_\rho(z)=\theta\left(\frac{|z-\bar{y}|^2}{\rho^2}\right)\;.
$$
Then, for any $\rho>0$ such that $B(\bar{y},\rho)\subset\subset {\cal O}$, we have
$$
\|\xi_\rho-\psi\|_\infty\leq C (\eta\rho+(M_1+M_2)\rho^2)\;,\; 
\|D\xi_\rho-D\psi\|_\infty\leq C(\eta+ (M_1+M_2)\rho)
$$
and 
$$
\|D^2\xi_\rho\|_\infty\leq \frac{C}{\rho}(\eta+ (M_1+M_2)\rho)
$$
for some constant $C=C({{k}})>0$, where we have set $\eta=|D\phi(\bar{y})-D\psi(\bar{y})|$, 
$M_1=\|D^2\phi\|_\infty$ and $M_2=\|D^2\psi\|_\infty$.
\end{Lemma}

\begin{Remark} \  \rm \\ 
1. The norms $\|\cdot\|_\infty$ are 
of course taken on $B(\bar{y}, \rho)$, since $\xi_\rho=\psi$ outside. \\
2. The key point is that $\xi_\rho$ coincides with $\phi$ in a small
neighbourhood of $\bar{y}$, but is not too far from $\psi$ in the full set ${\cal O}$
provided $\rho$ and $\eta$ are small.
\end{Remark}

We are now ready to prove the proposition. \\

\noindent {\bf Proof of Proposition \ref{turlututu}.} 
From Lemma \ref{pouet1} the function $\phi_n$ defined by
$$
\phi_n(z)= \langle Dd_K(x_n),z-y_n\rangle +\frac12 \langle (D^2d_K(x_n)-\frac{1}{n} I_{{k}})(z-y_n), z-y_n\rangle
$$
has a strict local maximum at $y_n$ on $K$.
Since $y_n$ is the unique projection of $x_n$ onto $K$, 
the map $\psi_n(z)=d_K(x_n)-|z-x_n|$ has a global strict maximum on $K$ at $y_n$. Hence
$$
\zeta_n(z)= \phi_n(z)\theta_n(z)+\psi_n(z)(1-\theta_n(z))
\qquad \mbox{\rm where} \ \theta_n(z)=\theta\left(\frac{|z-y_n|^2}{\rho_n^2}\right),
$$
has a global strict maximum at $y_n$ on $K$, provided we choose $\rho_n>0,$ $\rho_n\to 0$ and
$n$ large enough. 
Since $\phi_n$ is globally smooth with uniformly bounded second order derivative, and since
$\psi_n$ is smooth with uniformly bounded second order derivative ouside of the ball $B(x, d_K(x)/2)$,
and since finally $\psi_n(y_n)=\phi_n(y_n)=0$ and $D\psi_n(y_n)=D\psi_n(y_n)=p_n$, Lemma \ref{pouet2} states
that $\zeta_n$ and $D\zeta_n$ uniformly converge to the function $z\to d_K(x)-|z-x|$ and its derivative 
respectively, in the set $\R^{{k}}\backslash B(x, d_K(x)/2)$, and that $\|D^2\zeta_n\|_\infty$ is uniformly 
bounded in $\R^{{k}}\backslash B(x, d_K(x)/2)$.

Let us now denote  by ${\bf d}_{\Sigma}$ the signed distance to $\Sigma$ (see (\ref{def-dist-signee}) 
for a definition).
Since $\partial \Sigma$ is a ${\cal C}^{1,1}$ manifold, we can find some 
open neighbourhood ${\cal O}$ of $\partial \Sigma$ such that ${\bf d}_{\Sigma}$ 
is ${\cal C}^{1,1}$ in ${\cal O}$, with $\|D^2{\bf d}_{\Sigma}\|_\infty$ bounded in ${\cal O}$. 
For $z\in\R^k,$ we define
$$
d^n (z) = {\bf d}_{\Sigma} (z) - \beta_n |y-z|^2 \ \ \ {\rm where} \ \beta_n >0, \ \beta_n \to 0.
$$
Notice that $d^n$ is ${\cal C}^{1,1}$ in ${\cal O},$ that $d^n$ and its derivative converge
locally uniformly to ${\bf d}_{\Sigma}$ whereas $\|D^2 d^n\|_\infty$ is
bounded in ${\cal O}.$ The advantage of introducing $d^n$ is that $\{d^n \leq 0 \}$
is still externally tangent to $K$ at $y$ with $\partial K\cap \partial \{d^n \leq 0 \}=\{y\}$
(instead, $\Sigma$ can touch $K$ at many points). We claim that, if we choose $\beta_n =2 |y-y_n|^{1/3},$
then, at least for $n$ large,
\begin{equation}\label{alphan}
d^n (z)\leq d^n (y_n) \ \ \ {\rm for \ any } \ z\in K\backslash B(y_n, \beta_n /2).
\end{equation}
Indeed, for $z\in K\backslash B(y_n, \beta_n /2),$
\begin{eqnarray*}
d^n (z)-d^n (y_n) & \leq &  -\beta_n |y-z|^2 - {\bf d}_{\Sigma}(y_n) + \beta_n |y-y_n|^2 \\
& \leq &
-\beta_n (\beta_n / 2 -|y-y_n|)^2 + |y-y_n| +\beta_n |y-y_n|^2 \\
& \leq &
-|y-y_n|(1-4|y-y_n|^\frac23),
\end{eqnarray*}
which is nonpositive for large $n$ since $y_n\to y.$

We introduce the maps
$$
\xi_n(z)= \zeta_n(z)\tilde{\theta}_n(z)+ [d^n (z)-d^n (y_n)](1-\tilde{\theta}_n(z))\quad 
{\rm where }\; \tilde{\theta}_n(z)=\theta\left(\frac{|z-y_n|^2}{\sigma_n^2}\right)
$$
and
$$
\sigma_n=\max\{ |p_n-Dd^n (y_n)|, \beta_n/ \sqrt{2}\}\;.
$$
Let us notice that $\sigma_n\to0$ because $y_n\to y$ and $p_n$ and
$Dd^n (y_n) = D{\bf d}_{\Sigma}(y_n) -2\beta_n (y_n-y)$ converge both to
$p= D{\bf d}_{\Sigma}(y)$ since $D{\bf d}_{\Sigma}$ 
is continuous at $y$. We now use Lemma \ref{pouet2}, with
$\zeta_n(y_n)=[d^n (y_n)-d^n (y_n)]=0$ and $\eta_n:=|D\zeta_n(y_n)-Dd^n(y_n)|=
|p_n-D{\bf d}_{\Sigma}(y_n) +2\beta_n (y_n -y)|\to 0$. 
It states that $\xi_n$ and $D\xi_n$ uniformly converge to ${\bf d}_{\Sigma}$ and $D{\bf d}_{\Sigma}$, 
respectively. Moreover, since $\eta_n\leq \sigma_n,$ the second
order derivative of $\xi_n$ is uniformly bounded. 

Let us finally prove that 
the set $\Sigma_n=\{\xi_n<0\}$ satisfies our requirements. What we already proved on $\xi_n$ shows that
$\Sigma_n$ converges to $\Sigma$ (in the ${\cal C}^{1,{\rm b}}$ sense). From its construction, 
$\Sigma_n$ is smooth in a neighborhood of $y_n$, with normal at the point $y_n$ equal to $p_n$ 
and curvature equal to the restriction to
$(p_n)^\perp$ of $-(X_n-\frac{1}{n}I_{{k}})$. 

It remains to check that $\Sigma_n$ is externally tangent to
$K$ at $y_n$. It suffices to prove that $\xi_n(z)\leq 0$ for any $z\in K$, because
$\xi_n(y_n)=0$. Let $z\in K$. If $|z-y_n|\leq \sigma_n/\sqrt{2}$, then $\xi_n(z)=\zeta_n(z)\leq 0$ from the 
construction of $\zeta_n$.
If $|z-y_n|>\sigma_n/\sqrt{2}$, then $|z-y_n|>\beta_n/2$ and thus, from (\ref{alphan}), $d^n(z)\leq d^n(y_n)$. 
Since moreover $\zeta_n$ has a global maximum on $K$ at $y_n$, we finally have
$$
\begin{array}{rl}
\xi_n(z)\; = & \zeta_n(z)\tilde{\theta}_n(z)+ [d^n (z)-d^n (y_n)](1-\tilde{\theta}_n(z))\\
\leq & \zeta_n(y_n)\tilde{\theta}_n(z)+ [d^n (y_n)-d^n (y_n)](1-\tilde{\theta}_n(z))=0\;.
\end{array}
$$
In conclusion we have proved that $\xi_n$ has a global maximum on $K$ at $y_n$, and the proof is complete.
\QED

\subsection{Proof of Theorem \ref{interposition}}
\label{preuve-thm-interposition}

The first part of the theorem is an immediate consequence of Proposition \ref{troistrois} :
we set $\Sigma_1 =\Sigma -(t,x)+(\bar{s}_1 ,\bar{y}_1)$ and  $\Sigma_2 =\Sigma -(t,x)+(\bar{s}_2 ,\bar{y}_2)$
where $\Sigma$ and $(t,x)$ are given by Proposition \ref{troistrois} and we check that
$\Sigma_1$ and $\Sigma_2$ enjoy the desired properties.

Without loss of generality, we assume that $\delta \in (0,\bar{e}/(2+\bar{e}))$.
Let us introduce, for all $(\tau_1,z_1,\tau_2,z_2)\in\R^+\times\R^N\times \R^+\times \R^N,$
$$
f(\tau_1,z_1,\tau_2,z_2)= d^\sigma_{ C_1}(\tau_1,z_1)+d^\sigma_{\widehat{ C_2}}(\tau_2,z_2)
+\frac{1}{2\delta}|(\tau_1,z_1)-(\tau_2,z_2)|_\sigma^2\;.
$$
Then $d^\sigma_{ C_1}$, $d^\sigma_{\widehat{ C_2}}$ and $f$ are semi-concave functions 
in $(\R^+\times\R^N)^2\backslash ( C_1\cup \widehat{ C_2})^2$. We claim that 
$f$ has a minimum at $(\bar{\tau}_1,\bar{z}_1,\bar{\tau}_2,\bar{z}_2)$, where 
$$
(\bar{\tau}_1,\bar{z}_1)=\frac12(1+\delta) (\bar{s}_1,\bar{y}_1)+\frac12(1-\delta)(\bar{s}_2,\bar{y}_2)
$$
and
$$
(\bar{\tau}_2,\bar{z}_2)=\frac12(1-\delta) (\bar{s}_1,\bar{y}_1)+\frac12(1+\delta)(\bar{s}_2,\bar{y}_2).
$$
Indeed, on the one hand, an easy computation shows that $f(\bar{\tau}_1,\bar{z}_1,\bar{\tau}_2,\bar{z}_2)
=\bar{e}-\delta /2.$ On the other hand,
for all $(\tau_1,z_1,\tau_2,z_2)\in\R^+\times\R^N\times \R^+\times \R^N,$ we have
\begin{eqnarray*}
& & f(\tau_1,z_1,\tau_2,z_2) \\ 
& = & d^\sigma_{ C_1}( {\tau}_1,{z}_1) + d^\sigma_{\widehat{ C_2}}({\tau}_2,{z}_2)
+\frac{1}{2\delta} |(\tau_1,z_1)-(\tau_2,z_2)|_\sigma^2 \\
& = & |(\tau_1,z_1)-(t_1,x_1)|_\sigma +|(\tau_2,z_2)-(t_2,x_2)|_\sigma
+ \frac{1}{2\delta} |(\tau_1,z_1)-(\tau_2,z_2)|_\sigma^2 \\
& \geq & |(t_1,x_1)-(t_2,x_2)|_\sigma - |(\tau_1,z_1)-(\tau_2,z_2)|_\sigma
+ \frac{1}{2\delta} |(\tau_1,z_1)-(\tau_2,z_2)|_\sigma^2,
\end{eqnarray*}
where $(t_1,x_1)\in C_1$ and $(t_2,x_2)\in \widehat{ C_2}$ are the points for which
the distances $d^\sigma_{ C_1}( {\tau}_1,{z}_1)$ and $d^\sigma_{\widehat{ C_2}}({\tau}_2,{z}_2)$
are achieved. It follows that $f(\tau_1,z_1,\tau_2,z_2)\geq \bar{e}-\delta /2$
since $|(t_1,x_1)-(t_2,x_2)|_\sigma\geq \bar{e}$ and since, for all $r\geq 0,$ $-r+r^2/(2\delta)\geq -\delta/2.$
Finally, $(\bar{\tau}_1,\bar{z}_1,\bar{\tau}_2,\bar{z}_2)$ is a minimum for $f.$

Since the semi-concave function $f$ has a minimum at $(\bar{\tau}_1,\bar{z}_1,\bar{\tau}_2,\bar{z}_2)$,
Jensen maximum principle \cite{jensen88} (see also \cite{fs93}) states that one can find a sequence of points 
$(\bar{\tau}_{1,n},\bar{z}_{1,n},\bar{\tau}_{2,n},\bar{z}_{2,n})$ which converges to
$(\bar{\tau}_1,\bar{z}_1,\bar{\tau}_2,\bar{z}_2)$
and a nonnegative symmetric matrix $\bar{A}\in {\cal{S}}_{2N+2}$ such that the functions
$d^\sigma_{ C_1}$, $d^\sigma_{\widehat{ C_2}}$ and $f$ are twice differentiable 
at $(\bar{\tau}_{1,n},\bar{z}_{1,n})$, $(\bar{\tau}_{2,n},\bar{z}_{2,n})$
and $(\bar{\tau}_{1,n},\bar{z}_{1,n},\bar{\tau}_{2,n},\bar{z}_{2,n})$ respectively and such that
\begin{equation}\label{DfD2f}
Df(\bar{\tau}_{1,n},\bar{z}_{1,n},\bar{\tau}_{2,n},\bar{z}_{2,n})\to 0 \quad 
{\rm and }\quad D^2f(\bar{\tau}_{1,n},\bar{z}_{1,n},\bar{\tau}_{2,n},\bar{z}_{2,n}) \to \bar{A}\geq 0.
\end{equation}
In particular, since $Dd^\sigma_{ C_1}(\bar{\tau}_{1,n},\bar{z}_{1,n})$
and $Dd^\sigma_{\widehat{ C_2}}(\bar{\tau}_{2,n},\bar{z}_{2,n})$ exist, the projections of 
$(\bar{\tau}_{1,n},\bar{z}_{1,n})$ onto 
$ C_1$ and $\widehat{ C_2}$ respectively are unique, and equal to some $(\bar{s}_{1,n},\bar{y}_{1,n})$ 
and $(\bar{s}_{2,n},\bar{y}_{2,n}).$ Note that $(\bar{s}_1,\bar{y}_1)$ is the unique projection onto $ C_1$ of 
$(\bar{t},\bar{x})$, and therefore $(\bar{s}_{1,n},\bar{y}_{1,n})$ converges to $(\bar{s}_1,\bar{y}_1)$. 
For the same reason, $(\bar{s}_{2,n},\bar{y}_{2,n})$ converges to $(\bar{s}_2,\bar{y}_2)$. 
Since $d^\sigma_{ C_1}$ and $d^\sigma_{\widehat{ C_2}}$ are semi-concave 
functions, the matrices $D^2d^\sigma_{ C_1}(\bar{\tau}_{1,n},\bar{z}_{1,n})$ and 
$D^2d^\sigma_{\widehat{ C_2}}(\bar{\tau}_{2,n},\bar{z}_{2,n})$ are bounded from above:
namely (see for instance \cite[Proposition 22.2]{cs04})
\begin{equation}\label{EstiAbove1}
D^2_{xx}d^\sigma_{ C_1}(\bar{\tau}_{1,n},\bar{z}_{1,n})\leq \frac{1}{d^\sigma_{ C_1}(\bar{\tau}_{1,n},\bar{z}_{1,n})}
I_N
\end{equation}
and
\begin{equation}\label{EstiAbove2}
D^2_{xx}d^\sigma_{\widehat{ C_2}}(\bar{\tau}_{2,n},\bar{z}_{2,n})\leq 
\frac{1}{d^\sigma_{\widehat{ C_2}}(\bar{\tau}_{2,n},\bar{z}_{2,n})} I_N\;.
\end{equation}
Using
(\ref{DfD2f}), we get
\begin{equation}\label{D2f}
\left(\begin{array}{cc}
D_{xx}^2d^\sigma_{ C_1}(\bar{\tau}_{1,n},\bar{z}_{1,n})+\frac{1}{\delta} I_N & -\frac{1}{\delta} I_N\\
-\frac{1}{\delta} I_N& D_{xx}^2d^\sigma_{\widehat{ C_2}}(\bar{\tau}_{2,n},\bar{z}_{2,n})+\frac{1}{\delta} I_N
\end{array}\right)
\to A\geq 0, 
\end{equation}
where the matrix $A\in {\cal{S}}_{2N}$ is the restriction to $\R^N\times \R^N$ of $\bar{A}.$
In particular,  if we set $A=\left(\begin{array}{cc}A_1& A_3\\A_3&A_2\end{array}\right)$, then
$$
D_{xx}^2d^\sigma_{ C_1}(\bar{\tau}_{1,n},\bar{z}_{1,n})+
D_{xx}^2d^\sigma_{\widehat{ C_2}}(\bar{\tau}_{2,n},\bar{z}_{2,n})
+\frac{2}{\delta} I_N
\to A_1+A_2 \geq 0
$$
and therefore 
$D_{xx}^2d^\sigma_{ C_1}(\bar{\tau}_{1,n},\bar{z}_{1,n})$ and 
$D_{xx}^2d^\sigma_{\widehat{ C_2}}(\bar{\tau}_{2,n},\bar{z}_{2,n})$ are in fact bounded.
So, after relabelling all the sequences, we can assume that 
the restriction of $-D_{xx}^2d^\sigma_{ C_1}(\bar{\tau}_{1,n},\bar{z}_{1,n})$
to $(D_{x} d^\sigma_{ C_1}(\bar{\tau}_{1,n},\bar{z}_{1,n}))^\perp$
converges to some matrix $X_{1}$ while the restriction of
$D_{xx}^2d^\sigma_{\widehat{ C_2}}(\bar{\tau}_{2,n},\bar{z}_{2,n})$
to $(D_{x} d^\sigma_{\widehat{ C_2}}(\bar{\tau}_{2,n},\bar{z}_{2,n}))^\perp$
converges to some $X_{2}.$
Note that, from  (\ref{D2f}) we have 
$$
\left(\begin{array}{cc}
-X_1+\frac{1}{\delta} I_{N-1} & -\frac{1}{\delta} I_{N-1}\\
-\frac{1}{\delta} I_{N-1} & X_2+\frac{1}{\delta} I_{N-1}
\end{array}\right)\geq 0\;.
$$
Moreover, since 
$$
d^\sigma_{ C_1}(\bar{\tau}_{1,n},\bar{z}_{1,n})\to d^\sigma_{ C_1}(\bar{\tau}_{1},\bar{z}_{1})=\frac{\bar{e}(1-\delta)}{2}
$$
and
$$
d^\sigma_{\widehat{ C_2}}(\bar{\tau}_{2,n},\bar{z}_{2,n})\to d^\sigma_{\widehat{ C_2}}(\bar{\tau}_{2},\bar{z}_{2})
=\frac{\bar{e}(1-\delta)}{2}\;,
$$
we get from (\ref{EstiAbove1}, \ref{EstiAbove2})
$$
\left(\begin{array}{cc}
-X_1 & 0\\
0& X_2
\end{array}\right)\leq 
\frac{2}{\bar{e}(1-\delta)} I_{2(N-1)}\leq 
\frac{1}{\delta} I_{2(N-1)}
$$
because $\delta<\bar{e}/(\bar{e}+2).$ So we have proved (\ref{AA}). 

We now apply Proposition \ref{turlututu} to the sets $ C_1$ and $\Sigma_1$. 
Assumption (\ref{assumA}) holds since the set $\Sigma_1$ is externally tangent to 
$ C_1$ at $(\bar{s}_1,\bar{y}_1)$. Moreover, the point 
$(\bar{s}_1,\bar{y}_1)$ is the unique projection of the point $(\bar{\tau}_1,\bar{z}_1)$ onto $ C_1$. 
The points $(\bar{\tau}_{1,n},\bar{z}_{1,n})$
are points of twice differentiability of $d^\sigma_{ C_1}$, converge to $(\bar{\tau}_1,\bar{z}_1)$ and have
a unique projection $(\bar{s}_{1,n},\bar{y}_{1,n})$ onto $ C_1$. Therefore we can find a sequence of sets
$\Sigma_{1,n}$ with ${\cal C}^{1,1}$ boundary, such that $\Sigma_{1,n}$
is externally tangent to $ C_1$ at $(\bar{s}_{1,n},\bar{y}_{1,n})$ and has a ${\cal C}^2$ boundary 
near $(\bar{s}_{1,n},\bar{y}_{1,n})$. Note that, since $\Sigma_1$ is a ${\cal C}^{1,1}$ regular tube and since the sets 
$\Sigma_{1,n}$ converge to $\Sigma_1$ in the ${\cal C}^{1,{\rm b}}$ sense, $\Sigma_{1,n}$ are also 
${\cal C}^{1,1}$ regular tubes provided $n$ is sufficiently large.

From Proposition \ref{turlututu}(ii), we also have that the curvature matrix 
of $\Sigma_{1,n}$ at $(\bar{s}_{1,n},\bar{y}_{1,n})$ is equal to the restriction of 
$-(D^2d^\sigma_{ C_1}(\bar{\tau}_{1,n},\bar{z}_{1,n})+\frac{1}{n}I_{N+1})$ to the tangent
space of $\Sigma_{1,n}$ at $(\bar{s}_{1,n},\bar{y}_{1,n})$. 
Let us denote by $X_{1,n}$ the restriction of this curvature matrix to $\R^N$.
We notice that $X_{1,n}$ converges to $X_1$.

In the same way, applying Proposition \ref{turlututu} to the complementary of the tube 
$\Sigma_2$ which is externally tangent to $\widehat{C_2}$ at $(\bar{s}_2,\bar{y}_2),$
we can find a sequence of ${\cal C}^{1,1}$ tubes $\Sigma_{2,n}$ which are
externally tangent to $\widehat{C_2}$ at some points $(\bar{s}_{2,n},\bar{y}_{2,n})$,
such that $\Sigma_{2,n}$ are of class ${\cal C}^2$ near $(\bar{s}_{2,n},\bar{y}_{2,n})$,
and such that the curvature matrix $X_{2,n}$ to $\Sigma_{2,n}(\bar{s}_{2,n})$ 
at $\bar{y}_{2,n}$ converges to $X_2$. 
\QED

\section{The inclusion principle}

\subsection{Statement of the main theorem. Existence, uniqueness and stability}

 \begin{Theorem}[Inclusion principle]\label{coconuts}
Let $0\leq  \lambda_1<\lambda_2$ be fixed, 
${\cal K}_1$  be a subsolution of the FPP (\ref{fpp}) with speed $h_{\lambda_1}$
on the time interval 
$[0,T)$ for some $T>0$
and ${\cal K}_2$ be a supersolution on $[0,T)$ with speed $h_{\lambda_2}$.
If
$$
\overline{{\cal K}_1}(0)\cap \widehat{{\cal K}_2}(0)=\emptyset\;,
$$
then, for all $t\in[0,T),$
$$
\overline{{\cal K}_1}(t)\cap \widehat{{\cal K}_2}(t)=\emptyset\;.
$$
\end{Theorem}
Before proving Theorem \ref{coconuts}, we recall some applications of such inclusion
principle, omiting the proofs which are easy adaptations of those of 
\cite{cr04} and  \cite{cardaliaguet01}.  

Concerning the existence of solutions, we have the following

\begin{Proposition}\label{sol}
For any initial position $K_0$, with $S\subset {\rm int}(K_0)$ and $K_0$
bounded, there is (at least) one solution to the FPP (\ref{fpp})
for $h_\lambda$. 

Moreover, there is a largest solution
and a smallest solution to this problem.
The largest solution has a closed graph while the smallest solution
has an open graph in $\R^+\times\R^N$. The largest solution contains all the subsolutions
of the FPP (\ref{fpp}) with initial condition $K_0$, while the smallest solution is contained
in any supersolution.
\end{Proposition}

In general, one cannot expect to have a unique solution, i.e., the closure of the minimal
solution is not necessarily equal to the maximal solution. However, uniqueness is generic:

\begin{Proposition}
Let $(K_\lambda)_{\lambda>0}$ be a strictly increasing family of initial positions 
(i.e., $K_{\lambda'}\subset\subset K_\lambda$ for $0<\lambda'<\lambda$) such that
$K_\lambda\in {\cal{D}}$ for all $\lambda >0.$
Then the solution of 
the FPP (\ref{fpp}) for $h_\lambda$ with initial
position $K_\lambda$ is unique but for a countable subset of the $\lambda$'s.
\end{Proposition}

Stability of solutions is expressed by means of Kuratowski upperlimit of sets. Let us recall
that,
if $(A_n)_{n\in\N}$ is a sequence of sets in $\R^M,$ then the Kuratowski upperlimit
${A^* = \mathop{\rm Lim\, sup}_{n} A_n}$ is the subset of all
accumulation points of somes sequences of points in $(A_n)_{n\in\N},$ namely
\begin{eqnarray} \label{def-kura}
A^* & := &  \{ z\in \R^M : \exists (n_k)_{k\in\N} \ {\rm increasing \ sequence \ of \ integers}, 
\exists (z_k)_{k\in\N},  \nonumber \\
& & \hspace*{4.5cm} 
z_k\in A_{n_k} \ {\rm and} \ z= \mathop{\rm lim}_{k} y_k \}.
\end{eqnarray}
We define $A_*$ as the complementary of $\mathop{\rm Lim\, sup}_{n} \widehat{A_n}.$

\begin{Proposition}\label{stabilo} If ${\cal K}_n$ is a sequence of subsolutions for $h_\lambda$, locally uniformly
bounded with respect to $t$, then the Kuratowski upperlimit ${\cal K}^*= \mathop{\rm Lim\, sup}_{n} {\cal K}_n$ 
is also a subsolution for $h_\lambda$.

In a similar way, if ${\cal K}_n$ is a sequence of supersolutions for $h_\lambda$, locally uniformly
bounded with respect to $t$, then ${\cal K}_*$ is also a supersolution for $h_\lambda$.
\end{Proposition}

\subsection{Proof of Theorem \ref{coconuts}}
 
\begin{figure}[h]
\begin{center}
\epsfig{file=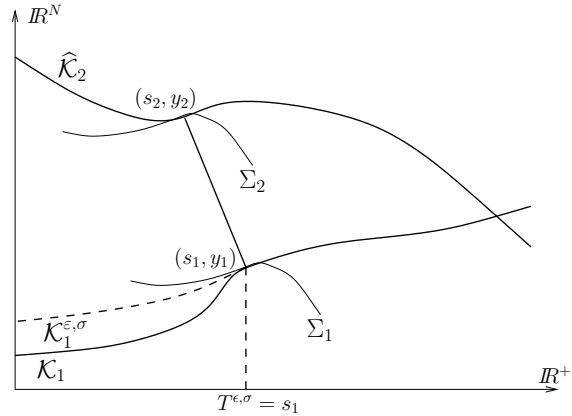, width=7.5cm} 
\vspace*{-0.7cm}
\end{center}
\caption{\label{dessin-preuv}
\textsl{Illustration of the proof of Theorem \ref{coconuts}.}}
\end{figure}

Without loss of generality, we assume that ${\cal K}_1$ has a  closed graph.
We argue by contradiction, assuming there exists $0\leq T^* < T$ such that
\begin{eqnarray} \label{hyp-absurde}
{\cal K}_1(T^*)\cap \widehat{{\cal K}_2}(T^*)\not=\emptyset.
\end{eqnarray}
For $\sigma>0$ and $\epsilon>0$, we consider 
$$
\rho_\sigma(t)=\min_{x\in{\cal K}_1(t)}d^\sigma_{\widehat{{\cal K}_2}}(t,x)
$$
and we set
$$
T^{\epsilon,\sigma}=\inf\{t\geq0\; : \; \rho_\sigma(t)\leq \epsilon e^{-t}\;\}\;.
$$
Recall that the notations $|\cdot|_\sigma$ and $d^\sigma_{\widehat{{\cal K}_2}}$
were introduced at the beginning of subsection \ref{InterpoTube}.
Let $r>0$ sufficiently small such that $S+rB$ has a ${\cal C}^2$ boundary and
$S+rB\subset\subset  {\cal K}_1(t)$ and $S+rB\subset\subset {\cal K}_2(t)$
for all $t\in [0,T]$. We also fix $R>0$ sufficiently large such that 
$$
\sup_{(t,x)\in {\cal K}_1, \; t\leq T}|x|+\sup_{(t,x)\in {\cal K}_2, \; t\leq T}|x|\leq R-r\;.
$$ 
We denote by $\theta$ 
the constant defined in Proposition \ref{regul} for $R$ and $r$.
Let us recall that $\theta>1/r$ and that, for any 
compact set $K$  with ${\cal C}^{1,1}$
boundary such that $S_r \subset {\rm int}(K)$ and $K\subset B(0,R-r )$, for any
$v\in \R^N$ with $|v|<1/\theta$ and any $x\in \partial K$, we have
\begin{equation}\label{IneqTranslat}
\bar{h}(x+v,K+v)\geq (1-\theta|v|)^2\bar{h}(x,K)\;.
\end{equation}
We refer the reader to Figure \ref{dessin-preuv} for an illustration of the
proof.

\begin{Lemma}\label{unun}  We can choose $\epsilon>0$ and $\sigma>0$
sufficiently small so that \\
(i) \ $\displaystyle{
\lambda_2(1-\theta \epsilon)^2\;>\lambda_1}$ and $\epsilon\sigma<1$ , \\ 
(ii) \ $T^{\epsilon,\sigma}>\epsilon ,$ \\
(iii) \ $\rho_\sigma(T^{\epsilon,\sigma}) = \epsilon e^{-T^{\epsilon,\sigma}},$ \\
(iv) for any $y_1\in\K1 (T^{\epsilon,\sigma})$ and any $(s_2,y_2)\in\HK2$ such that
\begin{equation}\label{Defy1y2}
\rho_\sigma(T^{\epsilon,\sigma}) =
\min_{y\in{\cal K}_1(T^{\epsilon,\sigma})}d^\sigma_{\widehat{{\cal K}_2}}(T^{\epsilon,\sigma},y)
= |(T^{\epsilon,\sigma},y_1)-(s_2,y_2)|_\sigma \, ,
\end{equation}
we have $y_1\neq y_2$ and $s_2>0$.
\end{Lemma}
The proof of the lemma is postponed to Section \ref{preuve-lemmes}.

From now on we fix $\epsilon$ and $\sigma$ as in Lemma \ref{unun} and we set $s_1=T^{\epsilon,\sigma}$.
Let us also set
\begin{eqnarray*}
\Kes &= &\{(s,y)\in\R^+\times\R^N : \\
& & \hspace*{0.3cm} \exists (\tau,z)\in \K1\ {\rm with}\ \tau\leq s_1\ {\rm and } \ 
|(s,y)-(\tau,z)|_\sigma\leq \epsilon (e^{-\tau}-e^{-s_1})\}.
\end{eqnarray*}
Recall that, for any two subsets $A_1$ and $A_2$ of $\R^{N+1}$, we define the minimal
distance between $A_1$ and $A_2$ by
$$
e(A_1,A_2)=\inf_{(t_1,x_1)\in A_1, \; (t_2,x_2)\in A_2} |(t_2,x_2)-(t_1,x_1)|_\sigma\;.
$$
\begin{Lemma} \label{deuxdeux} The set $\Kes$ is closed, with $\Kes(s_1)={\cal K}_1(s_1)$ and
$$
e(\Kes, \HK2)= \epsilon e^{-s_1}.
$$
Moreover there exist  $y_1\in\Kes (s_1)$ and $(s_2,y_2)\in \HK2$, such that
\begin{eqnarray} \label{hypdeuxdeux}
|(s_1,y_1)-(s_2,y_2)|_\sigma=e (\Kes, \HK2),
\end{eqnarray}
and, for such points $y_1$ and $(s_2,y_2)$, we have $y_1\neq y_2$.
\end{Lemma}
The proof is postponed.
Lemma \ref{deuxdeux} is a kind of refinement of Lemma \ref{unun}.(iv). 
Note that the proof of Lemma \ref{unun} and Lemma \ref{deuxdeux}
only use the fact that $\overline{{\cal K}_1}(0)\cap \widehat{{\cal K}_2}(0)=\emptyset$ and that
$\overline{{\cal K}_1}$ and $\widehat{{\cal K}_2}$ are left lower semicontinuous.

Next we give an estimate of the normal velocity of the tube
${\cal K}_1$ in terms of the normal velocity of $\Kes$:

\begin{Lemma}\label{TgModifTube} Assume that a ${\cal C}^{1,1}$ tube $\Sigma$ is externally tangent to $\Kes$
at some point $(\bar{s},\bar{y})\in \partial \Kes$. Let $(\bar{\tau},\bar{z})\in{\cal K}_1$ be such that
$$
\bar{\tau}\leq s_1\quad{\rm and}\quad |(\bar{s},\bar{y})-(\bar{\tau},\bar{z})|_\sigma
\leq \epsilon (e^{-\bar{\tau}}-e^{-s_1}).
$$
Then there is a ${\cal C}^{1,1}$ tube $\widetilde{\Sigma}$ which is externally tangent to the tube
${\cal K}_1\cap ([0,\bar{\tau}]\times \R^N)$ at $(\bar{\tau},\bar{z})$ and such that
$$
V^{\widetilde{\Sigma}}_{(\bar{\tau},\bar{z})}\geq V^\Sigma_{(\bar{s},\bar{y})}+\epsilon e^{-s_1}
\quad{\rm and}\quad 
\widetilde{\Sigma}(\bar{\tau})=\Sigma(\bar{s})+(\bar{z}-\bar{y}).
$$
If furthermore $\Sigma$ is of class ${\cal C}^2$ in a neighbourhood of $(\bar{s},\bar{y})$,
then $\widetilde{\Sigma}$ is also ${\cal C}^2$ in a neighbourhood of $(\bar{\tau},\bar{z})$.
\end{Lemma}
We postpone the proof. We are now ready to use the interposition Theorem \ref{interposition} that we apply to 
$C_1:=\Kes$ and $C_2:=\HK2$. Note that condition (\ref{Transversality}) holds thanks to 
Lemma \ref{deuxdeux}. Let us fix $(s_1,y_1)\in \Kes$ and $(s_2,y_2)\in \HK2$ with
$$|(s_1,y_1)-(s_2,y_2)|=e(\Kes,\HK2)\;.$$
From Theorem \ref{interposition} we know that 
there exists a regular tube $\Sigma_1$ with ${\cal C}^{1,1}$ boundary such that
$\Sigma_1$ is externally tangent to $\Kes$ at $(s_1,y_1)$
and $\Sigma_2:=\Sigma_1 +(s_2,y_2)-(s_1,y_1)$ is internally tangent to $\HK2$ at $(s_2,y_2)$
(see Figure \ref{dessin-preuv}). 

Futhermore, there exist ${\cal C}^{1,1}$ regular tubes
 $\Sigma_{1,n}$ and $\Sigma_{2,n}$ converging
respectively to $\Sigma_1$ and $\Sigma_2$ in the ${\cal C}^{1,{\rm b}}$ sense, there exist
$(s_{1,n}, y_{1,n})\in \Kes$ and $(s_{2,n}, y_{2,n})\in \HK2$ converging respectively to 
$(s_1,y_1)$ and $(s_2,y_2)$,  and there exist $(N-1)\times (N-1)$ matrices $X_1$, $X_2$, such that
\begin{itemize}
\item[(i)] $\Sigma_{1,n}$ is externally tangent to $\Kes$ at $(s_{1,n}, y_{1,n})$
and $\Sigma_{2,n}$ is internally tangent to $\HK2$ at $(s_{2,n},y_{2,n})$,

\item[(ii)] For $i=1,2$, $\Sigma_{i,n}$ is of class ${\cal C}^2$ in a neighbourhood 
of $(s_{i,n},y_{i,n})$ with
\begin{equation}\label{nunu1}
\lim_n\nu^{\Sigma_{1,n}(s_{1,n})}_{y_{1,n}}=\lim_n \nu^{\Sigma_{2,n}(s_{2,n})}_{y_{2,n}}= 
\nu^{\Sigma_1(s_1)}_{y_1}=\nu^{\Sigma_2(s_2)}_{y_2}, 
\end{equation} 
\begin{equation}\label{VV1}
\lim_n V^{\Sigma_{1,n}}_{(s_{1,n},y_{1,n})}=\lim_n V^{\Sigma_{2,n}}_{(s_{2,n},y_{2,n})}= 
V^{\Sigma_1}_{(s_1,y_1)}=V^{\Sigma_2}_{(s_2,y_2)}, 
\end{equation} 
and, for $i=1,2$, 
\begin{equation}\label{AA1}
\lim_n H^{\Sigma_{i,n}}_{s_{i,n},y_{i,n}}\to X_i,  \; 
{\rm with} \; X_1-X_2\leq 0\;.
\end{equation}
\end{itemize}

Since ${\cal K}_2$ is a supersolution for $h_{\lambda_2}$ and $\Sigma_{2,n}$ is internally tangent 
to ${\cal K}_2$ at $(s_{2,n},y_{2,n})$, we have
$$
V^{\Sigma_{2,n}}_{(s_{2,n},y_{2,n})}\geq F(\nu^{\Sigma_{2,n}(s_{2,n})}_{y_{2,n}},
H^{\Sigma_{2,n}(s_{2,n})}_{y_{2,n}}) +\lambda_2 \bar{h}(y_{2,n}, \Sigma_{2,n}(s_{2,n})).
$$
Letting $n\to+\infty$ gives, from the continuity property of $\bar{h}$, from (\ref{nunu1}), 
(\ref{VV1}) and (\ref{AA1})
\begin{equation}\label{Ineq1}
V^{\Sigma_2}_{(s_2,y_2)}\geq F(\nu^{\Sigma_2(s_2)}_{y_2},
X_2) +\lambda_2 \bar{h}(y_2, \Sigma_2(s_2))
\end{equation}

We now establish a symmetric inequality for $\Sigma_1$.
For any $n$, let $(\tau_{1,n},{z}_{1,n})\in {\cal K}_1$ be such that
$$
\tau_{1,n}\leq s_1\quad{\rm and}\quad 
|(s_{1,n},y_{1,n})-(\tau_{1,n},{z}_{1,n})|_\sigma\leq \epsilon (e^{-\tau_{1,n}}-e^{-s_1}).
$$
Note that $(\tau_{1,n}, {z}_{1,n})\to (s_1,y_1)$ because $(s_{1,n},y_{1,n})\to(s_1,y_1)$ and 
$\Kes$ has a closed graph with $\Kes(s_1)={\cal K}(s_1)$ (see Lemma \ref{deuxdeux}).
In particular, $\tau_{1,n}>0$ for large $n.$
Since $\Sigma_{1,n}$ is externally tangent to $\Kes$ at $(s_{1,n}, y_{1,n})$, Lemma
\ref{TgModifTube} states that one can find a ${\cal C}^{1,1}$ tube $\widetilde{\Sigma}_{1,n}$ 
externally tangent to ${\cal K}_1\cap([0, \tau_{1,n}]\times \R^N)$ at 
$(\tau_{1,n},{z}_{1,n})$ with 
$$
V^{\widetilde{\Sigma}_{1,n}}_{(\tau_{1,n},{z}_{1,n})}\geq 
V^{\Sigma_{1,n}}_{(s_{1,n},y_{1,n})}+\epsilon e^{-s_1}
\quad{\rm and}\quad 
\widetilde{\Sigma}(\tau_{1,n})=\Sigma(s_{1,n})+({z}_{1,n}-y_{1,n}).
$$
Moreover $\widetilde{\Sigma}_{1,n}$ is also ${\cal C}^2$ in a neighbourhood of $(\tau_{1,n},{z}_{1,n})$.
Lemma \ref{coupe}, given below, states that the tube ${\cal K}_1\cap([0, \tau_{1,n}]\times \R^N)$ 
is still a subsolution to the evolution equation for $h_{\lambda_1}.$
Therefore,
$$
V^{\widetilde{\Sigma}_{1,n}}_{(\tau_{1,n},{z}_{1,n})}\leq 
F(\nu^{\widetilde{\Sigma}_{1,n}(\tau_{1,n})}_{{z}_{1,n}}, 
H^{\widetilde{\Sigma}_{1,n}(\tau_{1,n})}_{{z}_{1,n}})
+\lambda_1 \bar{h}({z}_{1,n}, \widetilde{\Sigma}_{1,n}(\tau_{1,n}))\;.
$$
Since $\Sigma_{1,n}(s_{1,n})=\widetilde{\Sigma}_{1,n}(\tau_{1,n})-({z}_{1,n}-y_{1,n})$, 
we have $\nu^{\Sigma_{1,n}(s_{1,n})}_{y_{1,n}}=\nu^{\widetilde{\Sigma}_{1,n}(\tau_{1,n})}_{{z}_{1,n}}$
and $H^{\Sigma_{1,n}(s_{1,n})}_{y_{1,n}}=
H^{\widetilde{\Sigma}_{1,n}(\tau_{1,n})}_{{z}_{1,n}}$. Therefore
$$
V^{\Sigma_{1,n}}_{(s_{1,n},y_{1,n})}+\epsilon e^{-s_1}\leq 
F(\nu^{\Sigma_{1,n}(s_{1,n})}_{y_{1,n}}, 
H^{\Sigma_{1,n}(s_{1,n})}_{y_{1,n}})
+\lambda_1 \bar{h}( y_{1,n}, \Sigma_{1,n}(s_{1,n})-({z}_{1,n}-y_{1,n}))\;.
$$
Letting $n\to+\infty$ we get
\begin{equation}\label{Ineq2}
V^{\Sigma_1}_{(s_1,y_1)}+\epsilon e^{-s_1}\leq 
F(\nu^{\Sigma_1(s_1)}_{y_1}, X_1)
+\lambda_1 \bar{h}( y_1, \Sigma_1(s_1))\;.
\end{equation}

The difference between (\ref{Ineq1}) and (\ref{Ineq2}) gives
\begin{eqnarray}\label{Contradiction}
0 & \geq & F(\nu^{\Sigma_2(s_2)}_{y_2},X_2)+\lambda_2\bar{h}(y_2, \Sigma_2(s_2))
\nonumber \\
& & -F(\nu^{\Sigma_1(s_1)}_{y_1},X_1)-\lambda_1\bar{h}(y_1,\Sigma_1(s_1))
 +\epsilon e^{-s_1} \nonumber \\
& \geq & (\lambda_2(1-\theta\epsilon)^2-\lambda_1)\bar{h}(y_1,\Sigma(\bar{t})-e_1{\bar{\nu}}_x)
+\epsilon e^{-s_1},
\end{eqnarray}
because, on the one hand,
$V^{\Sigma_1}_{(s_1,y_1)}=V^{\Sigma_2}_{(s_2,y_2)}$,
$\nu^{\Sigma_1(s_1)}_{y_1}=\nu^{\Sigma_2(s_2)}_{y_2}$ and
$X_1\leq X_2$ (from (\ref{nunu1}, \ref{VV1},\ref{AA1})) and $F$ is elliptic and because, on the another hand, 
from (\ref{IneqTranslat}) and the definition of $\Sigma_2$,
$$
\bar{h}(y_2, \Sigma_2(s_2))=
\bar{h}(y_2, \Sigma_1(s_1)+y_2-y_1)
\geq (1-\theta|y_2-y_1|)^2\,\bar{h}(y_1,\Sigma_1(s_1)),
$$
where $|y_2-y_1|\leq \epsilon e^{-s_1}\leq \epsilon$. Since $\bar{h}\geq0$
and $\lambda_2(1-\theta\epsilon)^2\geq \lambda_1$ from Lemma \ref{unun}, 
there is a contradiction in (\ref{Contradiction}) and the proof is complete.
\QED

\subsection{Proofs of Lemmas \ref{unun}, \ref{deuxdeux}, \ref{TgModifTube} and \ref{coupe}}
\label{preuve-lemmes}

\noindent
{\bf Proof of Lemma \ref{unun}.}
The proof of this lemma is close to the one of \cite[Proposition 4.2]{cr04}. We provide it for
reader's convenience. 

The assertion (i) is obvious. 

To prove (ii), we first note that, since 
${\cal K}_1(0)\cap \widehat{{\cal K}_2}(0)=\emptyset$ and ${\cal K}_1$ and $\widehat{{\cal K}_2}$
are closed, there exists $\tau >0$ such that ${\cal K}_1(t)\cap \widehat{{\cal K}_2}(s)=\emptyset$
for all $0\leq s,t\leq \tau.$ Therefore
\begin{eqnarray*}
\mu := {\rm inf} \{ |y-x|  : 
0\leq s,t\leq \tau, \ x\in {\cal K}_1(t), \ y\in \widehat{{\cal K}_2}(s) \} >0.
\end{eqnarray*}
Let $\epsilon < \mu$ and $\sigma< \tau/(2\epsilon)$. Then, for $s\in[0,\tau/2]$ and $y\in{\cal K}_1(s)$
and for any $(t,x)\in \HK2$, we have
$$
|(s,y)-(t,x)|_\sigma\geq \left\{\begin{array}{ll}
\mu & {\rm if }\; t\leq \tau,\\
\frac{\tau}{2\sigma} & {\rm if}\;  t\geq \tau.
\end{array}\right.
$$
Hence $d_{\HK2}(s,y)\geq \epsilon$, which proves that $T^{\epsilon,\sigma}\geq \tau/2 \geq \epsilon$.

We prove (iii). From the definition of $\rho_\sigma (T^{\epsilon,\sigma}),$ there exists
$t_n \downarrow T^{\epsilon,\sigma}$ with $\rho_\sigma (t_n) \leq \epsilon e^{-t_n}.$
Therefore there exists $y_n\in {\cal K}_1(t_n)$ such that $\rho_\sigma (t_n) = 
d^\sigma_{\widehat{{\cal K}_2}}(t_n,y_n)$ and, up to extract a subsequence, we can
assume that $y_n \to y\in {\cal K}_1(T^{\epsilon,\sigma})$ (since ${\cal K}_1$
is closed). It follows
\begin{eqnarray*}
\rho_\sigma (t_n) = d^\sigma_{\widehat{{\cal K}_2}}(t_n,y_n) \to 
d^\sigma_{\widehat{{\cal K}_2}}(T^{\epsilon,\sigma},y)\geq \rho_\sigma (T^{\epsilon,\sigma}).
\end{eqnarray*}
Thus, we obtain the inequality $\rho_\sigma (T^{\epsilon,\sigma})\leq \epsilon e^{-T^{\epsilon,\sigma}}$
(note that we prove by the way that $\rho_\sigma$ is a lower-semicontinuous function). It remains
to prove that the equality holds. If not, there exists $y\in  {\cal K}_1(T^{\epsilon,\sigma})$
such that $d^\sigma_{\widehat{{\cal K}_2}}(T^{\epsilon,\sigma},y)< \epsilon e^{-T^{\epsilon,\sigma}}.$
From the left lower-semicontinuity of the subsolution ${\cal K}_1,$ for all sequence
$t_n \uparrow (T^{\epsilon,\sigma})^-,$ there exists $y_n\in {\cal K}_1(t_n)$ which converges to
$y.$ It follows that 
\begin{eqnarray*}
\rho_\sigma (t_n)\leq d^\sigma_{\widehat{{\cal K}_2}}(t_n,y_n) < \epsilon e^{-t_n} \ {\rm at \ least \
for\ }  n \ {\rm large.} 
\end{eqnarray*}
We get a contradiction with the definition of $T^{\epsilon,\sigma}$ and conclude for the proof of (iii).

We turn to the proof of (iv). For this we fix $\epsilon>0$ as in the proof of (ii) and we 
note that $T^{\epsilon,\sigma}$ is noncreasing with respect to $\sigma$. Since $T^{\epsilon,\sigma}\leq T^*$, 
$\lim_{\sigma\to 0^+}T^{\epsilon,\sigma}$ exists and we denote it by $\bar{s}$. Note that
$\bar{s}\geq \epsilon$.

We now argue by contradiction, assuming that there is a sequence $\sigma_n\to0^+$,
$y_{1,n}\in{\cal K}_1(T^{\epsilon,\sigma_n})$, $(s_{2,n},y_{2,n})\in\HK2$ such that
$$
\rho_{\sigma_n}(T^{\epsilon,\sigma_n}) 
= |(T^{\epsilon,\sigma_n},y_{1,n})-(s_{2,n},y_{2,n})|_{\sigma_n}=\epsilon e^{-T^{\epsilon,\sigma_n}} \, ,
$$
with either $s_{2,n}=0$ or $y_{1,n}=y_{2,n}$. Since ${\cal K}_1$ and $\widehat{{\cal K}_2}$ have closed graphs, up to extract subsequences,
there exist $y_1\in K_1(\bar{s})$ and $(s_2,y_2)\in \widehat{{\cal K}_2}$ such that
$y_{1,n}\to y_1$ and $(s_{2,n},y_{2,n})\to(s_2,y_2)$. From the inequality
$$
\epsilon\geq \left(\frac{1}{\sigma_n}(T^{\epsilon,\sigma_n}-s_{2,n})^2+|y_{1,n}-y_{2,n}|^2\right)^\frac12 \, ,
$$
we deduce that $T^{\epsilon,\sigma_n}-s_{2,n}\to 0$. But $T^{\epsilon,\sigma_n}\to \bar{s}>\epsilon$,
and so $s_{2,n}>0$ for $n$ sufficiently large, which implies that $y_{1,n}=y_{2,n}$
and thus $y_1=y_2$. 

We now use the left lower-semicontinuity of ${\cal K}_1$ and $\HK2$: Let $t_n\to \bar{s}^-$.
Since 
$y_1=y_2\in{\cal K}_1(\bar{s})\cap \HK2(\bar{s})$, 
there exist $x_{1,n}\in{\cal K}_1(t_n)$ and
$x_{2,n}\in\HK2(t_n)$
which converge to $y_1=y_2$. Then
$$
\rho_{\sigma_n}(t_n)\leq 
d^{\sigma_n}_{\HK2}(t_n, x_{1,n})\leq
|(t_n,x_{1,n})-(t_n, x_{2,n})|_\sigma
=|x_{1,n}-x_{2,n}|< \epsilon e^{-t_n}
$$
as soon as $n$ is sufficiently large. This is in contradiction with the definition of $T^{\epsilon, \sigma_n}$.
\QED

\noindent
{\bf Proof of Lemma \ref{deuxdeux}.} 
The set $\Kes$ is closed because so is ${\cal K}_1$. Let $y\in \Kes(s_1)$. There is some $(\tau,z)\in{\cal K}_1$
such that $\tau\leq s_1$ and $|(\tau,z)-(s_1,y)|_\sigma\leq \epsilon (e^{-\tau}-e^{-s_1})$. Then
$$
\frac{1}{\sigma^2}(\tau-s_1)^2\leq |(\tau,z)-(s_1,y)|_\sigma^2\leq \epsilon^2 (e^{-\tau}-e^{-s_1})^2\leq 
\epsilon^2(\tau-s_1)^2\;.
$$
Since $\epsilon\sigma <1$ (from Lemma \ref{unun}), we have  $\tau=s_1$ and thus $y=z$. 
So we have proved that $\Kes(s_1)\subset {\cal K}_1(s_1)$.
The other inclusion being obvious, the equality holds.

Let $(s,y)\in \Kes.$ Then there exists $(\tau,z)\in \K1$ with $\tau\leq s_1$ and
$|(s,y)-(\tau,z)|_\sigma \leq \epsilon (e^{-\tau}-e^{-s_1}).$ From the definition of
$\rho_\sigma$ and $s_1,$ we have 
$d^\sigma_{\widehat{{\cal K}_2}} (\tau,z)\geq \rho_\sigma (\tau) \geq \epsilon e^{-\tau}.$
It follows
\begin{eqnarray} \label{blabla11}
d^\sigma_{\widehat{{\cal K}_2}} (s,y) \geq 
d^\sigma_{\widehat{{\cal K}_2}} (\tau,z) - |(s,y)-(\tau,z)|_\sigma \geq \epsilon e^{-s_1}.
\end{eqnarray}
Taking the infimum over $(s,y)\in \Kes,$ we get $e(\Kes, \HK2)\geq \epsilon e^{-s_1}.$

Let us prove the opposite inequality. From Lemma \ref{unun}, we can choose
$y\in \K1 (s_1)$ such that $d^\sigma_{\widehat{{\cal K}_2}} (s_1,y)= \epsilon e^{-s_1}.$
But $\K1 (s_1)\subset \Kes(s_1).$ Therefore $e(\Kes, \HK2)\leq d^\sigma_{\widehat{{\cal K}_2}} (s_1,y)
= \epsilon e^{-s_1}.$

To prove the second assertion, let $y_1\in{\cal K}_1(s_1)$ and $(s_2,y_2)\in \HK2$
be such that (\ref{Defy1y2}) in Lemma \ref{unun} holds. Then obviously (\ref{hypdeuxdeux}) also holds.
For such points, let $(\tau,z)\in \K1$ be such that $\tau\leq s_1$ and 
\begin{eqnarray} \label{blabla12}
|(s_1,y_1)-(\tau,z)|_\sigma \leq \epsilon (e^{-\tau}-e^{-s_1}). 
\end{eqnarray}
If $\tau<s_1$, then, by definition of $s_1,$ we have $\rho_\sigma (\tau) > \epsilon e^{-\tau}.$
From the computation (\ref{blabla11}) with $(s,y):=(s_1,y_1),$ it follows 
$d^\sigma_{\widehat{{\cal K}_2}} (s_1,y_1)> \epsilon e^{-\tau} > \epsilon e^{-s_1}$ 
which is a contradiction with (\ref{hypdeuxdeux}). Hence $\tau=s_1$ and therefore
(\ref{blabla12}) gives $(s_1,y_1)=(\tau,z)\in \K1.$ We conclude by Lemma \ref{unun}.
\QED

\noindent
{\bf Proof of Lemma \ref{TgModifTube}.} Let ${\bf d}^\sigma_\Sigma$ be the signed distance to $\partial \Sigma$:
\begin{eqnarray} \label{DistSign}
{\bf d}^\sigma_\Sigma(\tau,z)=\left\{\begin{array}{ll}
d_\Sigma^\sigma (\tau,z) & {\rm if }\quad (\tau,z)\notin \Sigma, \\
-d_{\partial \Sigma}^\sigma (\tau,z) & {\rm if }\quad (t,z)\in \Sigma.
\end{array}\right.
\end{eqnarray}
Since $\Sigma$ is of class ${\cal C}^{1,1}$, one can find $\eta>0$ such that ${\bf d}^\sigma_\Sigma$
is of class ${\cal C}^{1,1}$ in $\{|{\bf d}^\sigma_\Sigma|<\eta\}$, with $D_x{\bf d}^\sigma_\Sigma\neq0$
in this set.

Let us define $\widetilde{\Sigma}$ by 
$$
\widetilde{\Sigma}=\{(\tau,z)\in\R^+\times \R^N\; : \; 
{\bf d}_{\Sigma}^\sigma((\tau,z)+(\bar{s},\bar{y})-(\bar{\tau},\bar{z}))\leq
-\epsilon (e^{-\tau}-e^{-\bar{\tau}})\}.
$$
We first show that the tube $\widetilde{\Sigma}$
is externally tangent to ${\cal K}_1\cap([0,\bar{\tau}]\times \R^N)$ at $(\bar{\tau},\bar{z}).$
At first, if $(\tau,z)\in {\cal K}_1$ with $\tau \leq \bar{\tau},$ then by definition of $\Kes,$
we have $B_\sigma ((\tau,z), \epsilon (e^{-\tau}-e^{-s_1}))\subset \Kes,$
where $B_\sigma$ is the usual open ball related to the norm $|\cdot |_\sigma.$
In particular, since $\Kes\subset \Sigma$, we have
${\bf d}_{\Sigma}^\sigma(\tau,z)\leq - \epsilon (e^{-\tau}-e^{-s_1})$. Therefore
$$
\begin{array}{l}
{\bf d}_{\Sigma}^\sigma((\tau,z)+(\bar{s},\bar{y})-(\bar{\tau},\bar{z})) \\
\hspace{1cm} \leq 
{\bf d}_{\Sigma}^\sigma(\tau,z)+|(\bar{s},\bar{y})-(\bar{\tau},\bar{z})|_\sigma \\
\hspace{1cm} \leq  - \epsilon (e^{-\tau}-e^{-s_1})+\epsilon (e^{-\bar{\tau}}-e^{-s_1})\\
\hspace{1cm}\leq - \epsilon (e^{-\tau}-e^{-\bar{\tau}})\;.
\end{array}
$$
So we have proved that ${\cal K}_1\cap([0,\bar{\tau}]\times \R^N)\subset \widetilde{\Sigma}$.
Moreover, we obviously have that $(\bar{\tau},\bar{z})\in \partial \widetilde{\Sigma}$
because $(\bar{s},\bar{y})\in \partial\Sigma$.

Let us show that $\widetilde{\Sigma}$ is a regular tube in an open interval containing
$\bar{\tau}$. For this, recall that ${\rm bd}(\Sigma)$ is defined by (\ref{DefBd}). If $(\tau,z)$ belongs to 
${\rm bd}( \Sigma)$ with $\epsilon |e^{-\tau}-e^{-\bar{\tau}}|<\eta$, then 
$$
|{\bf d}_{\Sigma}^\sigma((\tau,z)+(\bar{s},\bar{y})-(\bar{\tau},\bar{z}))|\; \leq
\epsilon |e^{-\tau}-e^{-\bar{\tau}}|<\eta
$$
and so ${\bf d}_{\Sigma}^\sigma$ is differentiable in a neighbourhood of 
$(\tau,z)+(\bar{s},\bar{y})-(\bar{\tau},\bar{z})$
with $D_x{\bf d}_{\Sigma}^\sigma\neq0$. Thus $\widetilde{\Sigma}$ is as smooth as $\Sigma$.
Finally,
$$
V^{\widetilde{\Sigma}}_{(\bar{\tau},\bar{z})}=
-\frac{\frac{\partial {\bf d}_{\Sigma}^\sigma}{\partial t}(\bar{s},\bar{y})
-\epsilon e^{-\bar{\tau}}}{|D_x{\bf d}_{\Sigma}^\sigma(\bar{s},\bar{y})|}
\geq V^\Sigma_{(\bar{s},\bar{y})}+\epsilon e^{-s_1}\;,
$$
since $|D_x{\bf d}_{\Sigma}^\sigma(\bar{s},\bar{y})|\leq |D{\bf d}_{\Sigma}^\sigma(\bar{s},\bar{y})|=1$.
\QED

\begin{Lemma}\label{coupe}
If ${\cal K}_1$ is a subsolution to the evolution equation for $h_{\lambda_1}$, then, for any
$t>0$, ${\cal K}_1\cap([0,t]\times \R^N)$ is also a subsolution for $h_{\lambda_1}$.
\end{Lemma}

\noindent
{\bf Proof of Lemma \ref{coupe}.} Let us set $\widetilde{{\cal K}_1}={\cal K}_1\cap([0,t]\times \R^N)$.
It is clear that $\widetilde{{\cal K}_1}$ is a left lower-semicontinuous tube because so is
${\cal K}_1$. Let $\Sigma$ be a ${\cal C}^2$ tube defined on some open time-intervall $I$ and
which is external tangent to $\widetilde{{\cal K}_1}$
at some point $(s,y)$ with $s>0$. If $s<t$, then, assuming without loss of generality that 
$I\subset (0,t)$, $\Sigma$ is also externally tangent to ${\cal K}_1$
and thus  $V^{\Sigma}_{(s,y)}\leq h_{\lambda_1}(y, \Sigma(s))\;.$

We now suppose that $s=t$ and, without loss of generality, 
that $\partial \Sigma\cap \widetilde{{\cal K}_1}=\{(s,y)\}$.
Let ${\bf d}^\sigma_\Sigma$ be the signed distance to $\partial \Sigma$
defined by (\ref{DistSign}). 
Note that, since $\partial \Sigma\cap \widetilde{{\cal K}_1}=\{(s,y)\}$,
${\bf d}^\sigma_\Sigma$ has a strict maximum on ${\cal K}_1$ at $(s,y)$ (at
least on the interval $I$). For $\gamma>0$ we introduce the mapping
$$
\varphi_\gamma(\tau,z)={\bf d}^\sigma_\Sigma(\tau,z)+\gamma \log(t-\tau)
$$
of class ${\cal C}^2$ for $\tau\in I\cap (0,t)$ and when $|{\bf d}^\sigma_\Sigma(\tau,z)|$
is small.
From standard arguments (see for instance the proof of Lemma 4.2 of \cite{cardaliaguet00}),
$\varphi_\gamma$ has a maximum on ${\cal K}_1$ at some point 
$(s_\gamma, y_\gamma)\in{\cal K}_1$ which converge to $(s,y)$ as $\gamma\to0^+$
and such that $s_\gamma<t$.
Moreover, the set 
$$
\Sigma_\gamma=\{(\tau, z)\in (I\cap(0,t))\times\R^N: \varphi_\gamma(\tau,z)\leq 
\varphi_\gamma(s_\gamma,y_\gamma)\}
$$
is a ${\cal C}^2$ tube on some open interval $I_\gamma\subset I\cap (0,t)$
with $s_\gamma\in I_\gamma$, and $\Sigma_\gamma(s_\gamma)$ converges in the ${\cal C}^2$
sense to $\Sigma(s)$. Now we note that $\Sigma_\gamma$ is externally tangent to ${\cal K}_1$
at $(s_\gamma,y_\gamma)$ and thus
$$
V^{\Sigma_\gamma}_{(s_\gamma,y_\gamma)}\leq h_{\lambda_1}(y_\gamma, \Sigma_\gamma(s_\gamma))\;,
$$
with 
$$
V^{\Sigma_\gamma}_{(s_\gamma,y_\gamma)}=
-\frac{ \frac{\partial \varphi}{\partial t}(s_\gamma,y_\gamma)}{|D\varphi(s_\gamma,y_\gamma)|}\;
\geq\; V^{\Sigma}_{(s_\gamma,y_\gamma)}\;.
$$
Hence
$$
V^{\Sigma}_{(s_\gamma,y_\gamma)} \leq h_{\lambda_1}(y_\gamma, \Sigma_\gamma(s_\gamma))\;,
$$
which gives the desired inequality as $\gamma\to0^+$: $V^{\Sigma}_{(s,y)}\leq h_{\lambda_1}(y, \Sigma(s))\;.$
\QED

\section{Convergence to equilibria}

In this section we investigate the asymptotic behavior of the solutions of our
front propagation problem (\ref{fpp}). More precisely, we show, under suitable assumptions
on the source $S$ and on $F$, that the solution converges as $t\to+\infty$ to the (weak) 
solution of the free boundary value problem:
\begin{eqnarray}\label{pb-bernoulli}
{\rm Find \ a \ set} \ K\in {\cal{D}} \ {\rm such \ that} \ 
h_\lambda(x,K)=0 \ {\rm for \ all} \  x\in \partial K
\end{eqnarray}
(recall that $\cal{D}$ is defined by (\ref{DefD})).
This problem is a generalization of the Bernoulli exterior free boundary problem
we recalled in the introduction.

Let us first introduce a notion of weak solution:

\begin{Definition} We say that the set $K\in {\cal{D}}$ is a viscosity
subsolution (respectively supersolution, solution) of the 
free boundary problem (in short FBP) (\ref{pb-bernoulli}) for $h_\lambda$ if the 
constant tube ${\cal K}(t)=K$ for all $t\geq0$ is a subsolution 
(respectively supersolution, solution) of the FPP (\ref{fpp})
for $h_\lambda$.
\end{Definition}
\begin{Remark} {\rm There are many other definitions of weak solutions for such 
FBP: see for instance the survey paper \cite{fr97}. 
The one we introduce here is the more suitable to our purpose. The idea of using
sub- and supersolution in this framework comes back to Beurling \cite{beurling58}.}
\end{Remark}

In order to ensure that the FBP (\ref{pb-bernoulli}) has a solution, we assume
in the sequel the following:
\begin{equation}\label{hyp-bernoulli}
\forall \lambda>0, \; \exists R>0\;\mbox{\rm such that } \ \forall r\geq R, \, \forall x\in B(0,r), \
h_\lambda(x, B(0, r))< 0.
\end{equation}
This assumption states that $B(0,r)$ is a strict classical supersolution of the free 
boundary problem for $h_\lambda$ for $r$ sufficiently large. It 
is in particular fulfilled (i) when $F(\nu,A)=Tr(A)+F_1(\nu)$, where $F_1(\nu)\leq 0$, or
(ii) when $F=F(\nu)<0$ for any $\nu$ with $|\nu|=1$, because of the behavior
of $\bar{h}$ for large balls (see Lemma \ref{growth}).
Note also that the assumption implies that, for any ball $B$, 
$$ 
F(\nu_x^B, H_x^B)<0\qquad \forall x\in \partial B\;,
$$
since $F$ is elliptic and $\bar{h}\geq0$.

\begin{Proposition}\label{ExistenceFBP} We assume that (\ref{Hypfg}) and (\ref{hyp-bernoulli}) hold.
Then, for any $\lambda>0$,  there is a largest and a smallest 
solution of the FBP (\ref{pb-bernoulli})
for $h_\lambda$, the largest being closed and containing any subsolution for $h_\lambda$, while
the smallest is open and is contained in all the supersolutions.
\end{Proposition}

\noindent{\bf Proof of Proposition \ref{ExistenceFBP}.} The proof can be achieved by a
direct application of Perron's method.
Existence and bounds for sub- and supersolutions are ensured by 
assumption (\ref{hyp-bernoulli}) and by Lemma \ref{SurSousSol}
below.
\QED

\begin{Lemma}\label{SurSousSol} We assume that (\ref{Hypfg}) and (\ref{hyp-bernoulli}) hold.
Let $\lambda>0$ be fixed. Then there exist $\epsilon>0$ and $R>0$ such that,
if $K$ is a subsolution (respectively a supersolution) of the FBP (\ref{pb-bernoulli}) for $h_\lambda$,
then $K\subset B(0,R)$ (respectively $S_\epsilon \subset K$ where $S_\epsilon$ is defined by (\ref{Sgamma})).
\end{Lemma}

\noindent{\bf Proof of Lemma \ref{SurSousSol}.} Let $r$ be the smallest nonnegative 
 real such that $K\subset B(0,r)$. Then there is some
point $x\in \partial K\cap \partial (B(0,r))$. Using the constant tube $B(0,r)\times \R$ as a test-tube and the fact that
$K$ is a subsolution, we have $h_\lambda(x,B(0,r))\geq0$, which implies that $r<R$ where $R$ is given by (\ref{hyp-bernoulli}).
Therefore we have $K\subset B(0,R)$. 

The assertion for the supersolution can be proved in a similar way, by using assumption (\ref{Hypfg})
and Lemma \ref{barhgrand} saying that $\bar{h}(x, S_\epsilon)$ is large for $\epsilon>0$ small and
$x\in\partial S_\epsilon.$
\QED

Next we address the uniqueness problem. The main assumption for this is that $S$
is starshaped. We also suppose that 
 $F=F(\nu,A)$ satisfies the subhomogeneity condition:
\begin{equation}\label{hyp-F2}
F(\nu,\gamma A)\geq \gamma F(\nu,A)\qquad \forall \gamma\geq 1,
\end{equation}
and that the following compatibility condition between $F$ and $S$ holds:
\begin{equation}\label{hyp-F3}
F(\nu^S_x, H^S_x)<0\qquad \forall x\in \partial S\;.
\end{equation}
Assumption (\ref{hyp-F2}) is fulfilled for instance if (i) $F(\nu,A)=Tr(A)+F_1(\nu)$, where $F_1(\nu)\leq 0$, or if
(ii) $F=F(\nu)< 0$ for any $\nu$ with $|\nu|=1$, while assumption (\ref{hyp-F3}) is always
satisfied for $F$ as in (ii), and is satisfied for sets with negative mean curvature if $F$ is as in (i).

\begin{Theorem}\label{unicite}  Let us assume that $S$ is strictly
starshaped at $0$, with $0\in {\rm int}(S)$, $g\equiv1$ and that 
(\ref{Hypfg}), (\ref{hyp-bernoulli}), (\ref{hyp-F2}) and (\ref{hyp-F3}) hold.
Then, for any $\lambda>0$, the FBP  (\ref{pb-bernoulli}) 
for $h_\lambda$ has a unique solution
denoted $K_\lambda$. 
Moreover, $K_\lambda$ is starshaped at $0$ for any $\lambda>0$
and the map $\lambda\to \overline{K_\lambda}$ is continuous for the Hausdorff topology.
\end{Theorem}
\begin{Remark} \rm $\;$\\
1. Uniqueness of solution means that, if $K_1$ and $K_2$ are two solutions of the FBP
for $h_\lambda$, then $\overline{K_1}=\overline{K_2}$ and ${\rm int}(K_1)={\rm int}(K_2)$.\\
2. Such a uniqueness result is classical in the literature, see in
particular Beurling \cite{beurling58}, Tepper \cite{tepper75}
and the survey paper \cite{fr97}. 
\end{Remark}
In order to prove Theorem \ref{unicite} we need three Lemmas. 
The first one explains that the homothetic of a subsolution is
still a subsolution. The second one allows to compare sub and supersolutions of the FBP.
The last one shows that subsolutions of the FBP for $h_\lambda$ 
when $\lambda$ is small are necessarily close to~$S$.

\begin{Lemma}\label{homot} Assume that $S$ is starshaped at $0$, $g\equiv1$ and that 
(\ref{hyp-F2}) holds. If $K$ is a subsolution of the FBP (\ref{pb-bernoulli}) for $h_\lambda$,
then $\rho K$ is a subsolution of the FBP for $h_{\rho \lambda}$
for any $\rho\in(0,1)$ such that $S\subset\subset \rho K$. 
\end{Lemma}

\noindent{\bf Proof of Lemma \ref{homot}.} 
For sake of clarity, we do the proof in a formal way, by assuming that $K$ is smooth.
If not, it is enough to do the same computations for the test-surfaces.
We first notice that 
\begin{equation}\label{ineq-bar-h}
\bar{h}(x,\rho K)\geq \frac{1}{\rho^2}\bar{h}(\frac{x}{\rho}, K)\qquad
\forall \rho\in (0,1), \ \forall x\in \partial (\rho K)\;.
\end{equation}
Indeed, if $u$ is the solution to (\ref{pde}) with $K$ instead of $\Omega$, then $v(x)=u(x/\rho)$ is a 
subsolution of equation (\ref{pde}) 
with $\rho K$ instead of $\Omega$ (we use here the fact that $S$ is starshaped, that $g\equiv1$
and, thus, that $0\leq u\leq 1$).
Then
$$
\bar{h}(x, \rho K)\geq |Dv(x)|^2=\frac{1}{\rho^2}|Du(x/\rho)|^2=\frac{1}{\rho^2}\bar{h}(x/\rho,K)\qquad
\forall x\in \partial( \rho K)\;.
$$
Next we also notice that $\nu^{\rho K}_x=\nu^K_{x/\rho}$ and $H^{\rho K}_x
=\frac{1}{\rho}H^K_{x/\rho}.$
Hence, for any $x\in \partial (\rho K)$, we have
$$
\begin{array}{rl}
h_{\rho\lambda}(x, \rho K)\;  = & F(\nu^{\rho K}_x,H^{\rho K}_x)+\rho\lambda \bar{h}(x, \rho K)\\
\geq & F(\nu^K_{x/\rho},\frac{1}{\rho}H^K_{x/\rho})+\frac{\lambda}{\rho} \bar{h}(x/\rho, K)\\
\geq & (1/\rho) h_\lambda(x/\rho, K) \geq 0
\end{array}
$$
because $K$ is a subsolution for $h_\lambda$. Hence $\rho K$ is a subsolution for $h_{\rho\lambda}$.
\QED

\begin{Lemma}\label{separation} We assume that (\ref{hyp-bernoulli}) holds.
Let $0<\lambda<\Lambda$, $R>0$ and $\gamma>0$ be fixed. 
Then, there is a constant $\kappa>0$, such that for any 
$\lambda\leq \lambda_1<\lambda_2\leq \Lambda$, for any 
subsolution $K_1$ of the FBP (\ref{pb-bernoulli}) for $h_{\lambda_1}$ and any supersolution $K_2$
for $h_{\lambda_2}$ with 
$$
S_\gamma \; \subset  \subset \; K_1\subset\subset K_2\; \subset \subset \; B(0,R-\gamma)\;,
$$
we have 
$$
K_1+\kappa (\lambda_2-\lambda_1)B(0,1)\subset  K_2 \;,
$$
where the sum in the above inclusion denotes the Kuratowski sum between sets. 
\end{Lemma}

\noindent{\bf Proof of Lemma \ref{separation}.} Let $\theta>0$ be the constant given
by Lemma \ref{regul}. From the assumption $K_1\subset\subset K_2$, 
we have $\overline{K_1}\cap \overline{\R^N\backslash K_2}=\emptyset$ and 
we can find $y_1\in \overline{K_1}$ and $y_2\in \overline{\R^N\backslash K_2}$ such that
$$
|y_1-y_2|=\min_{z_1\in K_1, z_2\in \overline{\R^N\backslash K_2}}|z_1-z_2|\;.
$$
Without loss of generality we can assume that $|y_1-y_2|<1/\theta$, since otherwise
the result is obvious.

\begin{figure}[h]
\begin{center}
\epsfig{file=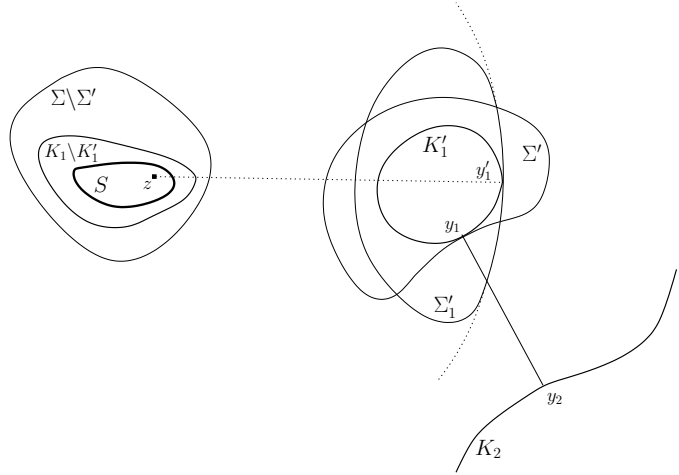, width=9cm} 
\vspace*{-0.7cm}
\end{center}
\caption{\label{dessin-lemme}
\textsl{Illustration of the proof of Lemma \ref{separation}.}}
\end{figure}

Using now the interposition and approximation results (see Proposition \ref{ilmanen} 
and Theorem \ref{app}), the fact that $K_1$
is a subsolution for $h_{\lambda_1}$ and $K_2$
a supersolution for $h_{\lambda_2}$, and proceeding as in the proof of Theorem
\ref{coconuts}, one can find an open set $\Sigma\subset \R^N$ with ${\cal C}^{1,1}$ boundary and 
$(N-1)\times (N-1)$ matrices $X_1\leq X_2$ such that 
\begin{equation}\label{111}
0\leq F(\nu_{y_1}^\Sigma, X_1)+\lambda_1 \bar{h}(y_1, \Sigma)
\end{equation}
and
\begin{equation}\label{222}
0\geq F(\nu_{y_2}^{\Sigma+y_2-y_1}, X_2)+\lambda_2 \bar{h}(y_2, \Sigma+y_2-y_1)\;.
\end{equation}
Since $\nu_{y_1}^\Sigma=\nu_{y_2}^{\Sigma+y_2-y_1}$ and $X_1\leq X_2$, we get
by subtracting (\ref{111}) to (\ref{222}) and using Lemma \ref{regul}
\begin{equation}\label{difference}
0\geq \left[\lambda_2\left(1-\theta|y_1-y_2|\right)^2-\lambda_1\right]\bar{h}(y_1,\Sigma)\;.
\end{equation}
In order to complete the proof, we have now to check that $\bar{h}(y_1,\Sigma)$
is positive. By Hopf maximum principle, we just have to show that the connected component
$\Sigma'$ of $\overline{\Sigma}$ which contains $y_1$ has a non empty intersection with the source $S$.
For this, we argue by contradiction, by assuming that $\Sigma'\cap S=\emptyset$ (see Figure
\ref{dessin-lemme} for an illustration). 
Let $K_1':=K_1\cap \Sigma'$. Note that 
\begin{equation}\label{separe}
e(K_1',\Sigma\backslash \Sigma')\geq e(\Sigma', \Sigma\backslash \Sigma')>0
\end{equation}
since $\Sigma$ is bounded with a smooth boundary. Let $z\in S$ and $y_1'$ be a point of maximum of the euclidean norm
$|\cdot-z|$ on $K_1'$. Note that $|y_1'-z|>0$.  The ball $B:=B(z, |y_1'-z|)$ is externally tangent
to $K_1'$ at $y_1'$. Thanks to (\ref{separe}), one can build an open set $\Sigma_1'$ with a ${\cal C}^2$
boundary, such that $K_1'\subset \overline{\Sigma_1'}$, $\overline{\Sigma_1'}\cap \overline{\Sigma\backslash \Sigma'}=\emptyset$,
$y_1'\in \partial \Sigma_1'$, and for which there
is a neighboourhood ${\cal O}$ of $y_1'$ with $\Sigma_1'\cap {\cal O}=B\cap {\cal O}$.
Note that $\overline{\Sigma_1'}\cap S=\emptyset$ since $S\subset\subset \Sigma\backslash \Sigma'$. 
Let us set $\Sigma_1=\Sigma_1'\cup (\Sigma\backslash \Sigma')$. Note that $\Sigma_1$ is externally tangent
to $K_1$ at $y_1'$. Moreover 
$$
h_\lambda (y_1', \Sigma_1)=
F(\nu_{y_1'}^{\Sigma_1},H_{y_1'}^{\Sigma_1})+\lambda\bar{h}(y_1',\Sigma_1)
= F(\nu_{y_1'}^B, H_{y_1'}^B)
$$
since $\bar{h}(y_1', \Sigma_1)=0.$ By (\ref{hyp-bernoulli}), for all $r\geq R$ and all
$x\in \partial B(0,r),$ we have $h_\lambda (x,B(0,r))<0.$ Therefore 
$F(\nu_x^{B(0,r)}, H_x^{B(0,r)})<0$ since $\bar{h}(x,B(0,r))\geq 0.$
By ellipticity, we have $F(\nu_{x'}^{B(0,r')}, H_{x'}^{B(0,r')})<0$ even for $r'\leq R$
and $|x'|=r'$ since $\nu_{x'}^{B(0,r')}= \nu_{rx'/|x'|}^{B(0,r)}.$ It follows
$$
h_\lambda (y_1', \Sigma_1)= F(\nu_{y_1'}^B, H_{y_1'}^B) <0
$$
which is a contradiction since $K_1$ is a subsolution. So $\bar{h}(y_1,\Sigma)~>~0.$

Then (\ref{difference}) leads to inequality
$$
|y_1-y_2|\geq \kappa(\lambda_2-\lambda_1)
$$
where $\kappa:=1/(2\theta\Lambda)$, which completes the proof.
\QED

\begin{Lemma}\label{CloseS} Under the assumptions of Theorem \ref{unicite},
for any $\epsilon>0$, there is $\lambda_0>0$ such that, 
for any subsolution $K$ of the FBP (\ref{pb-bernoulli}) for $h_\lambda$ with $\lambda \in(0,\lambda_0)$,
we have $K\subset S_\epsilon$ (see (\ref{Sgamma}) for a definition of $S_\epsilon$).
\end{Lemma}

\noindent{\bf Proof of Lemma \ref{CloseS}.} 
From assumption (\ref{hyp-F3}) and the regularity of the boundary of $S$, there is some $\alpha>0$
such that
$$
F(\nu^S_x, H^S_x)\leq -\alpha \qquad \forall x\in \partial S\;.
$$
Let us notice that a similar inequality also holds for $\rho S$, for $\rho\geq 1$, because
$$
F(\nu^{\rho S}_x, H^{\rho S}_x)=F(\nu^S_{x/\rho},\frac{1}{\rho}H^S_{x/\rho})
\leq \frac{1}{\rho} F(\nu^S_{x/\rho},H^S_{x/\rho})\leq -\frac{\alpha}{\rho}\;,
$$
thanks to assumption (\ref{hyp-F2}). 

Let us now fix $\epsilon>0$ and $\rho_0>1$ such that $\rho_0S\subset S_\epsilon$. 
Note that $S\subset\subset \rho_0S$ since $S$ is strictly starshaped.
We set $\lambda_0=\alpha/(\beta\rho_0)$, where $\beta=\sup_{x\in \partial (\rho_0S)}\bar{h}(x, \rho_0S)$.
Let $K$ be a subsolution for $h_\lambda$ with $\lambda\in(0,\lambda_0)$.
We denote by $\rho>1$ the smallest real such that $K\subset \rho S$. In order to prove that
$\rho\leq \rho_0$, we argue by contradiction and assume that $\rho>\rho_0$. Since 
$\rho S$ is externally tangent to $K$ at some point $x\in\partial K$ and $K$
is a subsolution, we have
$$
0\leq h_\lambda(x, \rho S) =F(\nu^{\rho S}_x, H^{\rho S}_x)+\lambda \bar{h}(x, \rho S)
\leq -\frac{\alpha}{\rho}+\lambda \bar{h}(x, \rho S)
$$
where, from inequality (\ref{ineq-bar-h}), 
$$
\bar{h}(x, \rho S)\leq \left(\frac{\rho_0}{\rho}\right)^2\bar{h}(\frac{\rho_0x}{\rho}, \rho_0S)\leq
\left(\frac{\rho_0}{\rho}\right)^2\beta\;.
$$
Hence $0\leq-\frac{\alpha}{\rho}+\lambda \left(\frac{\rho_0}{\rho}\right)^2\beta,$
which implies that $\rho\leq \lambda_0\rho_0^2\beta/\alpha=\rho_0$, a contradiction.
So we have proved that $\rho\leq \rho_0$. Therefore $K\subset \rho_0 S\subset S_\epsilon$.
\QED

\noindent{\bf Proof of Theorem \ref{unicite}.} Let us denote for any $\lambda>0$
by $K_\lambda$ the maximal solution of the BFP for $h_\lambda$. Note that
$\lambda\to K_\lambda$ is nondecreasing since $K_\lambda$ contains all the
subsolutions for $h_\lambda$.

We first check that $K_\lambda$ is starshaped at $0$. Indeed, from Lemma \ref{homot},
for any $\rho\in(0,1)$ sufficiently close to 1, the set $\rho K_\lambda$
is a subsolution for $h_{\rho\lambda}$, and thus for $h_\lambda$. Since $K_\lambda$
contains all the subsolutions, we have $\rho K_\lambda\subset K_\lambda$ for any 
$\rho\in(0,1)$ sufficiently close to 1, which implies that $K_\lambda$ is starshaped.

Next, we show that the map $\lambda\to K_\lambda$ is continuous
for the Hausdorff topology. From the stability property of solutions (Proposition \ref{stabilo}), 
the decreasing limit of the $K_{\lambda'}$
converges to $K_\lambda$ when $\lambda'\to \lambda^+$. Hence we only have to show
that ${\rm Lim}_{\lambda'\to\lambda^-}K_{\lambda'}$ equals $K_\lambda$,
where Lim denotes the Kuratowski limit (see (\ref{def-kura})).

Since, for any $\rho\in(0,1)$ sufficiently close to 1, 
the set $\rho K_\lambda$ is a subsolution for $h_{\rho\lambda}$, 
we have $\rho K_\lambda\subset K_{\rho\lambda}$.
Therefore
$$
K_\lambda=\mathop{\rm Lim}_{\rho\to 1^-} \rho K_\lambda\subset 
\mathop{\rm Lim}_{\lambda'\to \lambda^-}  K_{\lambda'} \subset K_\lambda.
$$
So we have checked that $\lambda\to K_\lambda$ is continuous.

Let us finally prove that, for any $\lambda>0$, $K_\lambda$ is the unique
solution of
for $h_\lambda$. Let $K$ be another solution. Note that $K\subset
K_\lambda$. 
From Lemma \ref{CloseS}, we can find 
some $\overline{\lambda}_1>0$ such that $K_{\overline{\lambda}_1}\subset\subset
K$ because $S\subset\subset K$.
Let us set 
$$
\overline{\lambda}=\sup\{ \lambda'\; |\; K_{\lambda'}\subset\subset K\}\;.
$$
We now use Lemma \ref{separation} with $r>0$ and $R$ such that 
$S_r\subset K_{\overline{\lambda}_1}$ and $K_\lambda\subset B(0,R).$
There is a constant $\kappa>0$ such that for any 
$\overline{\lambda}_1<\lambda'<\overline{\lambda}$, 
$$
K_{\lambda'}+\kappa(\lambda-\lambda')B\subset  K \;.
$$
The continuity of the map $\lambda'\to K_{\lambda'}$ then implies that 
$$
K_{\overline{\lambda}}+\kappa(\lambda-\overline{\lambda})B\subset  \overline{K}
\;.
$$
Therefore $\overline{\lambda}=\lambda$ since, otherwise, the continuity of 
$\lambda'\to K_{\lambda'}$ would also imply the existence of $\epsilon>0$
such that $K_{\overline{\lambda}+\epsilon}\subset\subset \overline{K}$, a contradiction
with the definition of $\overline{\lambda}$. Therefore $K_\lambda\subset\overline{K}.$

In order to prove that ${\rm int}(K)={\rm int}(K_\lambda)$,
we notice that ${\rm int}(K_\lambda)=\bigcup_{\lambda'<\lambda} K_{\lambda'}$,
because $K_{\lambda'}\subset\subset K_\lambda$ for $\lambda'<\lambda$,
and therefore the equality $\bar{\lambda}=\lambda$ implies that ${\rm int}(K_\lambda)\subset {\rm int}(K)$.
Since the converse inequality is obvious, the proof of Theorem \ref{unicite} is complete.
\QED

\begin{Corollary}[Asymptotic behavior]\label{AsBe} Under the assumptions of Theorem \ref{unicite},
if ${\cal K}$ is a solution of the FPP (\ref{fpp}) for $h_\lambda$, then ${\cal K}(t)$
converges, for the Hausdorff metric as $t\to +\infty,$  to the unique solution $K_\lambda$ 
of the FBP (\ref{pb-bernoulli}) for $h_\lambda$ while $\widehat{K}(t)$ converges to 
$\overline{\R^N\backslash K}$.
\end{Corollary}

\begin{Remark} {\rm \ \\
1. Note that the above result holds for any solution ${\cal K}(t)$ of the FPP (\ref{fpp})
with any initial position ${\cal{K}}(0) \in {\cal{D}}.$ \\
2. The proof of the asymptotic behavior which follows relies strongly on the uniqueness
of the solution of the limit problem (\ref{pb-bernoulli}).    }
\end{Remark}

\noindent{\bf Proof of Corollary \ref{AsBe}.} Let us fix $\lambda_1<\lambda<\lambda_2$.
From lemma \ref{barhgrand} and (\ref{hyp-bernoulli}), there are
$r>0$ and $R>0$ such that $S_r$ and $B(0,R)$ are respectively sub- and 
supersolution to the FBP for $h_{\lambda_1}$
and $h_{\lambda_2}$. We can also choose $r>0$ sufficiently small and $R>0$
sufficiently large so that $S_r\subset\subset {\cal
K}(0)\subset\subset B(0,R)$. The inclusion principle then states that 
$$
S_r\subset\subset{\cal K}(t)\subset\subset B(0,R)\qquad \forall
t\geq0\;.
$$

Let $K^*$ be the Kuratowski upperlimit of ${\cal K}(t)$  as $t\to+\infty$
(see (\ref{def-kura})). 
Note that $S_r\subset K^*\subset B(0,R)$ and that
the constant tube $\R\times K^*$ is actually the upperlimit of the solutions
${\cal K}(\cdot+\tau)$ as $\tau\to+\infty$.
From the stability of solutions  (see Proposition \ref{stabilo}),
the constant tube $\R\times K^*$ is a subsolution of the FPP for $h_\lambda$.
Hence, $K^*$ is a subsolution of the FBP (\ref{pb-bernoulli}) for $h_\lambda$ and
we have  $K^*\subset K_\lambda$. 

In the same way, if we set $L^*$ to be the Kuratowski upperlimit of
$\widehat{K}(t)$ as $t\to +\infty,$ then $\R^N\backslash L^*$ is a
supersolution to FBP (\ref{pb-bernoulli}) for $h_\lambda.$ Since
$K_\lambda$ is the unique solution for $h_\lambda$, $K_\lambda$ is also the smallest
solution, which implies that $K_\lambda\subset \overline{\R^N\backslash L^*}.$
Hence 
$$
K^*\subset K_\lambda\subset \overline{\R^N\backslash L^*}.
$$
Since $\R^N\backslash L^*\subset K^*$, the proof is complete.
\QED


\end{document}